\documentclass[aos,MSNbibl,nameyear,seceqn,dvips]{arximspdf}
\usepackage{accents}
\usepackage{subeqn}
\usepackage{graphicx}

\doi{10.1214/13-AOS1173} 
\volume{42}
\issue{1}
\pubyear{2014}
\firstpage{29}
\lastpage{63}

\makeatletter
\def\smash{}
\renewcommand{\mathring}[1]{\accentset{\circ}{#1}}

\newcommand{\rrvert}{\vert}
\newcommand{\llvert}{\vert}
\newcommand{\eqref}[1]{(\ref{#1})}
\makeatletter
\newcommand{\inlineoverset}[2]{%
{\mathop{#2}\limits^{#1}}}
\makeatother
\newproclaim{definition}{Definition}[section]
\newproclaim{example}{Example}[section]
\newtheorem{theorem}{Theorem}[section]
\newtheorem{proposition}{Proposition}[section]
\newtheorem{corollary}{Corollary}[section]
\newtheorem{lemma}{Lemma}[section]
\newproclaim{remark}{Remark}[section]
\makeatother

\begin{document}
\begin{frontmatter}

\title{Co-clustering separately exchangeable network~data\thanksref{T1}}
\runtitle{Co-clustering network data}
\thankstext{T1}{Supported in part by the US Army Research Office under
PECASE Award W911NF-09-1-0555 and MURI Award W911NF-11-1-0036; by the UK
EPSRC under Mathematical Sciences Established Career Fellowship
EP/K005413/1 and Institutional Sponsorship Award EP/K503459/1; by the
UK Royal Society under a Wolfson Research Merit Award; and by Marie
Curie FP7 Integration Grant PCIG12-GA-2012-334622 within the 7th
European Union Framework Program.}

\begin{aug}
\author[a]{\fnms{David} \snm{Choi}\ead[label=e1]{davidch@andrew.cmu.edu}}
\and
\author[b]{\fnms{Patrick J.} \snm{Wolfe}\corref{}\ead[label=e2]{p.wolfe@ucl.ac.uk}}

\runauthor{D. Choi and P. J. Wolfe}

\affiliation{Carnegie Mellon University and University College London}

\address[a]{Heinz College of Public Policy and Management\\
Carnegie Mellon University\\
5000 Forbes Ave, Hamburg Hall\\
Pittsburgh, Pennsylvania 15213-3890\\
USA\\
\printead{e1}}

\address[b]{Department of Statistical Science\\
University College London\\
Gower Street\\
London WC1E 6BT\\
United Kingdom\\
\printead{e2}}

\end{aug}

\received{\smonth{12} \syear{2012}}
\revised{\smonth{9} \syear{2013}}

%
\begin{abstract}
This article establishes the performance of stochastic blockmodels in
addressing the co-clustering problem of partitioning a binary array
into subsets, assuming only that the data are generated by a
nonparametric process satisfying the condition of separate
exchangeability.\vspace*{1pt} We provide oracle inequalities with rate of
convergence $\mathcal{O}_P (n^{-1/4} )$ corresponding to profile
likelihood maximization and mean-square error minimization, and show
that the blockmodel can be interpreted in this setting as an optimal
piecewise-constant approximation to the generative nonparametric model.
We also show for large sample sizes that the detection of co-clusters
in such data indicates with high probability the existence of
co-clusters of equal size and asymptotically equivalent connectivity in
the underlying generative process.
\end{abstract}

%


\end{frontmatter}

\section{Introduction}

Blockmodels are popular tools for network modeling that see wide and
rapidly growing use in analyzing social, economic and biological
systems; see \citet{zhao2011community} and \citet{fienberg2012brief}
for recent overviews. A blockmodel dictates that the probability of
connection between any two network nodes is determined only by their
respective block memberships, parameterized by a latent categorical
variable at each node.

Fitting a blockmodel to a binary network adjacency matrix yields a
clustering of network nodes, based on their shared proclivities for
forming connections. More generally, fitting a blockmodel to any binary
array involves partitioning it into blocks. In this way, blockmodels
represent a piecewise-constant approximation to a latent function that
generates network connection probabilities. This in turn can be viewed
as a histogram-like approximation to a nonparametric generative process
for binary arrays; fitting such models is termed co-clustering [\citet
{flynn2012consistent,rohe2012co}].

This article analyzes the performance of stochastic blockmodels for
co-clustering under model misspecification, assuming only an underlying
generative process that satisfies the condition of separate
exchangeability [\citet{diaconis2007graph}]. This significantly
generalizes known results for the blockmodel and its co-clustering
variant, which have been established only recently under the
requirement of correct model specification [\citet{bickel2009nonparametric,bickelmethod,rohe2010spectral,chatterjee2012matrix,choi2012stochastic,flynn2012consistent,rohe2012co,zhao2013consistency,fishkind2013consistent}].

We show that blockmodels for co-clustering satisfy consistency
properties and remain interpretable whenever separate exchangeability
holds. Exchangeability is a natural condition satisfied by many network
models: it characterizes permutation invariance, implying that the
ordering of nodes carries no information [\citet
{bickel2009nonparametric,hoff2009cmot}]. A blockmodel is an
exchangeable model in which the connection probabilities are piecewise
constant. Blockmodels also provide a simplified parametric
approximation in the more general nonparametric setting [\citet{bickelmethod}].

In addition to providing oracle inequalities for blockmodel
$M$-estimators corresponding to profile likelihood and least squares
optimizations, we show that it is possible to identify clusterings in
data---what practitioners term \emph{network communities}---even when
the actual generative process is far from a blockmodel. The main
statistical application of our results is to enable co-clustering under
model misspecification. Much effort has been devoted to the task of
community detection [\citet{newman2006modularity,fortunato2007resolution,zhao2011community, fienberg2012brief}], but the
drawing of inferential conclusions in this setting has been limited by
the need to assume a correctly specified model.

Our results imply that community detection can be understood as finding
a best piecewise-constant or simple function approximation to a
flexible nonparametric process. In settings where the underlying
generative process is not well understood and the specification of
models is thus premature, such an approach is a natural first step for
exploratory data analysis. This has been likened to the use of
histograms to characterize exchangeable data in nonnetwork settings
[\citet{bickel2009nonparametric}].

The article is organized as follows. In Section~\ref{sec:modelElicit},
we introduce our nonparametric setting and model. In Section~\ref{sec:oracle} we present oracle inequalities for co-clustering based on
blockmodel fitting. In Section~\ref{sec:consistency} we give our main
technical result, and discuss a concrete statistical application:
quantifying how the collection of co-clusterings of the data approaches
that of a generative nonparametric process. We prove our main result in
Section~\ref{subsec:sketch}, by combining a construction used to
establish a theory of graph
limits [\citeauthor{borgs2006graph}
(\citeyear{borgs2006graph,borgs2008convergent,borgs2007convergent})] with statistical learning
theory results on $U$-statistics [\citet{clemenccon2008ranking}]. In
Section~\ref{subsec:simulations} we illustrate our results via a
simulation study, and in Section~\ref{sec:discussion} we relate them to
other recent work. Appendices \ref{sec:proof.th2}--\ref
{sec:simulationLemmaPf} contain additional proofs and technical lemmas.

\section{Model elicitation}
\label{sec:modelElicit}

Recall that fitting a blockmodel to a binary array involves
partitioning it into blocks. Denote by $G = (V_1,V_2,E)$ a bipartite
graph with edge set $E$ and vertex sets $(V_1,V_2)$, where assignments
of vertices to $V_1$ or $V_2$ are known. For example, $V_1$ and $V_2$
might represent people and locations, with edge $(i,j)$ denoting that
person $i$ frequents location $j$. See \citet{flynn2012consistent} and
\citet{rohe2012co} for additional examples.

\subsection{Exchangeable graph models}

For a bipartite graph $G$ represented as a binary array $A$, the
appropriate notion of exchangeability is as follows.

%
\begin{definition}[({Separate exchangeability [\citet{diaconis2007graph}]})]
\label{def:exchangeable}
An array $\{A_{ij}\}_{i,j=1}^\infty$ of binary random variables is
separately exchangeable if
\[
P(A_{ij} = X_{ij}, 1\leq i,j \leq n) = P(A_{ij} =
X_{\Pi_1(i)\Pi_2(j)}, 1\leq i,j \leq n)
\]
for all $n = 1, 2, \ldots,$ all permutations $\Pi_1,\Pi_2$ of
for all $n = 1, 2, \ldots,$ all permutations $\Pi_1,\Pi_2$ of
$1,\ldots
,n$, and all $X \in\{0,1\}^{n \times n}$.
\end{definition}

If we identify a finite set of rows and columns of $A$ with the
adjacency matrix of an observed bipartite graph $G$, then it is clear
that the notion of separate exchangeability encompasses a broad class
of network models. Indeed, given a single observation of an unlabeled
graph, it is natural to consider the class of all models that are
invariant to permutation of its adjacency matrix; see \citet
{bickel2009nonparametric} and \citet{hoff2009cmot} for discussion.

The assumption of separate exchangeability is the only one we will
require for our results to hold. A representation of models in this
class will be given by the Aldous--Hoover theorem for separately
exchangeable binary arrays.

%
\begin{definition}[(Exchangeable array model)] \label{model:general}
Fix a measurable mapping $\omega\dvtx  [0,1]^3 \rightarrow[0,1]$. Then the
following model generates an exchangeable random bipartite graph
$G=(V_1,V_2,E)$ through its adjacency matrix $A$:
\begin{longlist}[(1)]
\item[(1)] generate $\alpha\sim\operatorname{Uniform}(0,1)$;
\item[(2)] fix $m = |V_1|$ and $n = |V_2|$, and generate each element
of $\xi=(\xi_1,\ldots,\xi_m)$ and $\zeta=(\zeta_1,\ldots,\zeta_n)$
$\inlineoverset{\mathrm{i.i.d.}}{\sim}
\operatorname{Uniform}(0,1)$;\vadjust{\goodbreak}
\item[(3)] for $i = 1,\ldots,m$, and $j=1,\ldots,n$, generate $A_{ij}
\inlineoverset{\mathrm{i.i.d.}}{\sim} \operatorname{Bernoulli} (\omega^\alpha
(\xi_i,\zeta_j) )$, where $\omega(x,y) \equiv\omega^\alpha(x,y)$
denotes the function $(x,y) \mapsto\omega(\alpha,x,y)$. If $A_{ij}=1$,
then connect vertices $i\in V_1$ and $j \in V_2$.
\end{longlist}
\end{definition}

The Aldous--Hoover theorem states that this representation is
sufficient to describe any separately exchangeable network distribution.

\begin{theorem}[{[\citet{diaconis2007graph}]}]\label{th:exchangeable graph}
Let $\{A_{ij}\}_{i,j=1}^\infty$ be a separately exchangeable binary
array. Then there exists some $\omega\dvtx  [0,1]^3 \rightarrow[0,1]$,
unique up to measure-preserving transformation, which generates $\{
A_{ij}\}_{i,j=1}^\infty$.
\end{theorem}

The interpretation of the exchangeable graph model of Definition~\ref
{model:general} is that each vertex has a latent parameter in $[0,1]$
($\xi_i$ for vertex $i$ in $V_1$, and $\zeta_j$ for vertex $j$ in
$V_2$) which determines its affinity for connecting to other vertices,
while $\alpha$ is a network-wide connectivity parameter
(nonidentifiable from a single network observation). Because $\xi$ and
$\zeta$ are latent, $\omega(x,y)$ itself is identifiable only up to
measure-preserving transformation, and is hence indistinguishable from
any mapping $(x,y) \mapsto\omega(\alpha,\pi_1(x),\pi_2(y))$ for which
$\pi_1,\pi_2$ are in the set $\mathcal{P}$ of measure-preserving
bijective maps of $[0,1]$ to itself.

\subsection{The stochastic co-blockmodel}

Many popular network models can be recognized as instances of
Definition~\ref{model:general}. For example, \citet{hoff2002latent,airoldi2008mixed} and
\citet{kim2012multiplicative} all present models
in which the resulting $\omega(\alpha,x,y)$ is constant in $\alpha$,
while \citet{miller2009nonparametric} require the full parameterization
$\omega(\alpha,x,y)$. The stochastic co-blockmodel specifies $\omega
(\alpha,x,y)$ constant in $\alpha$ and also piecewise-constant in $x$
and $y$, and thus can be viewed as a simple function approximation to
$\omega(x,y)$ in Definition~\ref{model:general}.

%
\begin{definition}[{(Stochastic co-blockmodel [\citet{rohe2012co}])}]\label
{model:co-block}
Fix integers $K_1,K_2 > 0$, a matrix $\theta\in[0,1]^{K_1 \times
K_2}$ and discrete probability measures $\mu$ and $\nu$ on $\{
1,\ldots,K_1\}$ and $\{1,\ldots,K_2\}$. Then the \emph{stochastic
co-block\-model} generates an exchangeable bipartite graph
$G=(V_1,V_2,E)$ through the matrix $A$ as follows:
\begin{longlist}[(1)]
\item[(1)] Fix $m = |V_1|$ and $n = |V_2|$, and generate $S =
(S(1),\ldots,S(m)) \inlineoverset{\mathrm{i.i.d.}}{\sim} \mu$ and $T =
(T(1),\ldots
,T(n)) \inlineoverset{\mathrm{i.i.d.}}{\sim} \nu$.
\item[(2)] For $i=1,\ldots,m$, and $j=1,\ldots,n$, generate $A_{ij}
\inlineoverset{\mathrm{i.i.d.}}{\sim} \operatorname{Bernoulli} (\theta_{S(i)
T(j)} )$. If $A_{ij} =1$, then connect vertices $i \in V_1$ and $j
\in V_2$.
\end{longlist}
Additionally, given co-blockmodel parameters $\phi\equiv(\mu,\nu
,\theta)$, define
\[
\omega_\phi(x,y) = \theta_{F^{-1}_\mu(x) F^{-1}_\nu(y)},\qquad x,y \in[0,1]
\]
as the mapping corresponding to Definition~\ref{model:general}, with
$F_\mu^{-1}(x) = \inf_z \{ F_\mu(z) \geq x \}$ the inverse distribution
function corresponding to a given distribution $\mu$.\vadjust{\goodbreak}
\end{definition}

Without loss of generality we assume $K_1 = K_2 = K$ in what follows,
noting that our results do not depend in any crucial way on this
assumption. Thus, a stochastic blockmodel's vertices in $V_1$ belong to
one of $K$ latent classes, as do those in~$V_2$. Vectors $S \in\{
1,\ldots,K\}^m$ and $T \in\{1,\ldots,K\}^n$ of categorical variables
specify these class memberships. The matrix $\theta\in[0,1]^{K \times
K}$ indexes the corresponding connection affinities between classes in
$V_1$ and $V_2$. Because $S$ and $T$ are latent, the stochastic
co-blockmodel is identifiable only up to a permutation of its class labels.

\section{Oracle inequalities for co-clustering}
\label{sec:oracle}

If we assume that the separately exchangeable data model of
Definition~\ref{model:general} is in force, then a natural first step
is to approximate $\omega(x,y)$ by way of some piecewise-constant
$\omega_\phi(x,y)$, according to the stochastic co-blockmodel of
Definition~\ref{model:co-block}. This approximation task is equivalent
to fixing $K$ and estimating $\phi= (\mu,\nu,\theta)$ by co-clustering
the entries of an observed adjacency matrix $A \in\{0,1\}^{m \times n}$.

\subsection{Sets of co-clustering parameters}

To accomplish this task, we consider $M$-estimators that involve an
optimization over the latent categorical variable vectors $S \in\{
1,\ldots,K\}^m$ and $T \in\{1,\ldots,K\}^n$. The resulting blockmodel
estimates will reside in a set $\Phi$ containing triples $(\mu,\nu
,\theta) \in\Omega_m \times\Omega_n \times[0,1]^{K \times K}$, where
we define $\Omega_m$ to be the set of all probability distributions
over $\{1,\ldots,K\}$ whose elements are integer multiples of $1/m$,
\[
\Omega_m = \Biggl\{ p \in \biggl\{ 0,\frac{1}{m},
\frac{2}{m},\ldots,1 \biggr\}^K \dvtx \sum
_{a=1}^K p_a = 1 \Biggr\},
\]
and likewise for $\Omega_n$. Note that $\Omega_m$ and $\Omega_n$ are
subsets of the standard $K-1$-simplex, chosen to contain all measures
$\mu$ and $\nu$ that can be obtained by empirically co-clustering the
elements of an $m \times n$-dimensional binary array. Thus, by
construction, any estimator $\hat{\phi}(A) = (\hat{\mu},\hat{\nu
},\hat
{\theta})$ based on an empirical co-clustering of an observed binary
array $A \in\{0,1\}^{m \times n}$ has codomain $\Phi$.

Given a specific $\mu$ and $\nu$, let $\mathcal{Q}_\mu^m$ denote the
set of all node-to-class assignment functions that partition the set $\{
1,\ldots,m\}$ into $K$ classes in a manner that respects the
proportions dictated by $\mu= (\mu_1, \ldots, \mu_K) \in\Omega_m$,
\[
\mathcal{Q}_\mu^m = \bigl\{v \in\{1,\ldots,K
\}^m\dvtx \bigl|v^{-1}(a)\bigr| = m \mu_a, a=1,\ldots,K
\bigr\},
\]
and likewise for $\mathcal{Q}_\nu^n$.

\subsection{Oracle inequalities}

We now establish that, for $L^2$ risk and Kullback--Leibler divergence,
there exist $M$-estimators that enable us to determine, with rate of
convergence $n^{-1/4}$, optimal piecewise-constant approximations of
the generative $\omega(x,y)$, up to quantization due to the
discreteness of $\Phi$.

\begin{theorem}[(Oracle inequalities for co-clustering)]\label{th:learning}
Let $A \in\{0,1\}^{m \times n}$ be a separately exchangeable array
generated by some $\omega$ in accordance\vadjust{\goodbreak} with Definition~\ref
{model:general}, and consider fitting a $K$-class stochastic
co-blockmodel parameterized by $\phi\equiv(\mu,\nu,\theta)$ to $A$.
Then as $n \rightarrow\infty$, with $K$ and $m/n$ fixed:

\begin{longlist}[(1)]

\item[(1)] For the least squares co-blockmodel $M$-estimator
%
\begin{equation}
\label{eq:MSfunctional} \hat{\phi} = \mathop{\operatorname{argmin}}_{\phi
\in\Phi} \Biggl\{ \min
_{S \in\mathcal
{Q}_\mu^m,T \in\mathcal{Q}_\nu^n} \frac{1}{mn} \sum_{i=1}^m
\sum_{j=1}^n \llvert
\theta_{S(i) T(j)} - A_{ij} \rrvert^2 \Biggr\}
\end{equation}
relative to the $L^2$ risk
\[
R_\omega(\phi) = \inf_{\pi_1,\pi_2 \in\mathcal{P}} \int_{[0,1]^2}
\bigl\llvert\omega \bigl(\pi_1(x),\pi_2(y) \bigr) -
\omega_\phi(x,y) \bigr\rrvert^2 \,dx \,dy,
\]
we have that
\[
R_\omega(\hat{\phi}) - \inf_{\phi\in\Phi} R_\omega(
\phi) = \mathcal{O}_P \bigl(n^{-1/4} \bigr);
\]

\item[(2)] Given any $\phi= (\mu,\nu,\theta)$, let $B(\phi) =
\max_{1\leq a,b,\leq K} | \log(\theta_{ab}/(1-\theta_{ab})) |$.
Consider the profile likelihood co-blockmodel $M$-estimator
%
\begin{eqnarray}
\label{eq:PLfunctional} &&\hat{\phi} = \mathop{\operatorname{argmax}}_{\phi
\in\Phi} \Biggl\{ \max
_{S \in
\mathcal{Q}_\mu^m, T \in\mathcal{Q}_\nu^n} \frac{1}{mn} \sum_{i=1}^m
\sum_{j=1}^n \bigl\{ A_{ij}
\log(\theta_{S(i)T(j)})
\nonumber
\\[-8pt]
\\[-8pt]
\nonumber
&&\hspace*{163pt}{}+ (1-A_{ij}) \log(1-\theta_{S(i) T(j)}) \bigr\} \Biggr\}
\end{eqnarray}
relative to
\begin{eqnarray*}
&&L_\omega(\phi) = \sup_{\pi_1,\pi_2 \in\mathcal{P}} \int_{[0,1]^2}
\bigl\{ \omega \bigl(\pi_1(x),\pi_2(y) \bigr) \log
\omega_\phi(x,y)\\
&&\hspace*{106pt}{}+ \bigl[1-\omega \bigl(\pi_1(x),\pi_2(y) \bigr) \bigr]
\log \bigl(1-\omega_\phi(x,y) \bigr) \bigr\} \,dx \,dy.
\end{eqnarray*}
If $\phi^* = \operatorname{argmax}_{\phi\in\Phi} L_\omega(\phi)$
exists, and $B(\phi^*)$ and $B(\hat{\phi})$ are finite, then
%
\begin{equation}
\label{eq:th:learning:1} \frac{\max_{\phi\in\Phi} L_\omega(\phi
) - L_\omega(\hat{\phi})}{B(\phi
^*) + B(\hat{\phi})} = \mathcal{O}_P
\bigl(n^{-1/4} \bigr).
\end{equation}
\end{longlist}
\end{theorem}

Theorem~\ref{th:learning} can be viewed as analyzing maximum likelihood
techniques in the context of model misspecification [\citet
{white1982maximum}], and is proved in Appendix \ref{sec:proof.th2}. It
establishes that minimization of the squared error between a fitted
co-blockmodel and an observed binary array according to \eqref
{eq:MSfunctional} serves as a proxy for approximation of $\omega$ by
$\omega_\phi$ in mean square, and that fitting a stochastic
co-blockmodel via profile likelihood according to \eqref
{eq:PLfunctional} is equivalent to minimizing the average
Kullback--Leibler divergence of the approximation $\omega_\phi(x,y)$
from the generative $\omega(x,y)$.

The existence of a limiting object $\omega(x,y)$ implies that we are in
the dense graph regime, with expected network degree values increasing
linearly as a function of~$m$ or $n$. Given a correctly specified
generative blockmodel, profile likelihood estimators are known to be
consistent even in the sparse graph setting of polynomial or
poly-logarithmic expected degree growth [\citet
{bickel2009nonparametric}]. In our setting, however, the generative
model is no longer necessarily a blockmodel; in this context, both
\citet{borgs2008convergent} and \citet{chatterjee2012matrix} leave open
the question of consistently estimating sparse network parameters,
while \citet{bickelmethod} give an identifiability result extending to
the sparse case. The simulation study reported in Section~\ref{subsec:simulations} below suggests that the behavior of blockmodel
estimators is qualitatively similar across at least some families of
dense and sparse models.

\subsection{\texorpdfstring{Additional remarks on Theorem \protect\ref{th:learning}}
{Additional remarks on Theorem 3.1}}

In essence, Theorem~\ref{th:learning} implies that the binary array $A$
yields information on its underlying generative $\omega(x,y)$ at a rate
of at least $n^{-1/4}$. While the necessary optimizations in \eqref
{eq:MSfunctional} and \eqref{eq:PLfunctional} are not currently known
to admit efficient exact algorithms, they strongly resemble existing
objective functions for community detection for which many authors have
reported good heuristics [\citet{newman2006modularity,fortunato2007resolution,zhao2011community}]. Furthermore,
polynomial-time spectral algorithms are known in certain settings to
find correct labelings under the assumption of a generative blockmodel
[\citet{rohe2010spectral,fishkind2013consistent}], suggesting that
efficient algorithms may exist when distinct clusterings or community
divisions are present in the data. In this vein, \citet
{chatterjee2012matrix} has recently proposed a universal thresholding
procedure based on the singular value decomposition.

\begin{remark}
We may replace the objective function of \eqref{eq:PLfunctional}\break  with
the full profile likelihood function $\max_{S \in\mathcal{Q}_\mu^m, T
\in\mathcal{Q}_\nu^n}
\{ \sum_{i=1}^m \log\mu_{\smash{S(i)}} +\break \sum_{j=1}^n \log\nu
_{\smash{T(j)}} +  \sum_{i=1}^m \sum_{j=1}^n \{ A_{ij} \log\theta
_{S(i)T(j)} + (1-A_{ij}) \log(1-\theta_{S(i) T(j)}) \} \}$. The
same rate of convergence can then be established with respect to the
corresponding term for $L_\omega(\phi)$, adapting the proofs in
Appendices \ref{sec:proof.th2} and \ref{sec:proof.th3}.
\end{remark}

\begin{remark}\label{rem:Bterms}
Assume $ \phi^* = \operatorname{argmax}_{\phi\in\Phi} L_\omega
(\phi)$
exists.\vspace*{1pt} Terms $B(\phi^*)$ and $B(\smash{\hat{\phi}})$ in \eqref
{eq:th:learning:1} show that elements of $\theta^*$ and $\smash{\hat
{\theta}}$ must not approach $0$ or $1$ too quickly as $n \rightarrow
\infty$; otherwise $L_\omega(\hat{\phi})$ can be much smaller than
$L_\omega(\phi^*)$.
\end{remark}

This is a natural consequence of the fact that the Kullback--Leibler
divergence of $\omega_\phi$ from $\omega$ is finite if and only if
$\omega$ is absolutely continuous with respect to $\omega_\phi$. To see
the implication, consider $\xi,\zeta$, and $A$ generated according to
Definition~\ref{model:general} with $\omega(x,y) = 1\{x\leq1/2\} 1\{
y\leq1/2\}$. Let $\mu_1 = m^{-1} \sum_{i=1}^m 1\{\xi_i \leq1/2\}$ and
$\nu_1 = n^{-1} \sum_{j=1}^n 1\{\zeta_i \leq1/2\}$. Then the
maximum-likelihood two-class blockmodel fit to~$A$ will yield $\omega
_{\smash{\hat{\phi}}}(x,y) = 1\{x \leq\mu_1\} 1\{y \leq\nu_1\}$,
and so $L_\omega(\smash{\hat{\phi}})$ diverges to $-\infty$ unless
$\mu
_1 = \nu_1 = 1/2$.

\section{Convergence of co-cluster estimates}
\label{sec:consistency}

We now give our main technical result and show its statistical
application in enabling us to interpret the convergence of co-cluster
estimates. The estimators of Theorem~\ref{th:learning} require
optimizations over the set of all possible co-clusterings of the data;
that is, over vectors $S$ and $T$ that map the observed vertices to
$1,\ldots,K$. Analogously, one may also envision an uncountable set of
co-clusterings of the generative model, which map the unit interval
$[0,1]$ to $1,\ldots,K$. We define these two sets of co-clusterings
more formally and then give a result showing in what sense they become
close with increasing $m$ and $n$, so that optimizing over co-clusters
of the data is asymptotically equivalent to optimizing over co-clusters
of the generative model. This result yields the rate of convergence
$\mathcal{O}_P(n^{-1/4})$ appearing in Theorem~\ref{th:learning}, and
also has a geometric interpretation that sheds light on the estimators
defined by \eqref{eq:MSfunctional} and \eqref{eq:PLfunctional}.

\subsection{\texorpdfstring{Relating co-clusterings of $A$ to those of $\omega$}
{Relating co-clusterings of $A$ to those of omega}}
\label{sec:notation}

Given a bipartite graph $G=(V_1,V_2,E)$ with adjacency matrix $A \in\{
0,1\}^{m \times n}$, recall that the latent class vectors $S \in\{
1,\ldots,K\}^m$ and $T \in\{1,\ldots,K\}^n$ respectively partition
$V_1$ and $V_2$ into $K$ subsets each. To relate an empirical
co-clustering of $A$ to a piecewise-constant approximation of some
$\omega$, we first define the matrix $A/ST \in[0,1]^{K \times K}$ to
index the proportion of edges spanning each of the $K^2$ subset pairs
defined by $S$ and $T$,
\[
(A/ST)_{ab} = \frac{1}{mn} \sum_{i \in S^{-1}(a)}
\sum_{j \in
T^{-1}(b)} A_{ij},\qquad a,b=1,\ldots,K.
\]
Second, we define mappings $\sigma,\tau\dvtx [0,1]\rightarrow\{1,\ldots
,K\}
$, which will play a role analogous to $S$ and $T$. Given some $\omega
\dvtx [0,1]^2 \rightarrow[0,1]$, this allows us to define a matrix $\omega
/\sigma\tau\in[0,1]^{K \times K}$ which encodes the mass of $\omega$
assigned to each of the $K^2$ subset pairs defined by $\sigma$ and
$\tau
$ as follows:
\[
(\omega/\sigma\tau)_{ab} = \int_{\sigma^{-1}(a) \times\tau
^{-1}(b)} \omega(x,y)
\,dx \,dy,\qquad a,b=1,\ldots,K.
\]

We will use the $K \times K$ matrices $A/ST$ and $\omega/\sigma\tau$
to index all possible co-clusterings that can be induced by
partitioning an observed binary array $A \in\{0,1\}^{m \times n}$ into
$K^2$ blocks. To link these sets of co-clusters, recall from
Section~\ref{sec:oracle} the sets $\mathcal{Q}_\mu^m$ and $\mathcal
{Q}_\nu^n$ of all node-to-class assignment functions that partition $\{
1,\ldots,m\}$ and $\{1,\ldots,n\}$ into $K$ classes in manners that
respect the proportions dictated by $\mu= (\mu_1, \ldots, \mu_K)
\in
\Omega_m$ and $\nu= (\nu_1, \ldots, \nu_K) \in\Omega_n$. Analogously,
we define $\mathcal{Q}_\mu$ (resp., $\mathcal{Q}_\nu$) to be the set of
partitions of $[0,1]$ into $K$ subsets whose cardinalities are of
proportions $\mu_1,\ldots,\mu_K$:
\begin{eqnarray*}
\mathcal{Q}_\mu&= \bigl\{\sigma\dvtx [0,1]\rightarrow\{1,\ldots,K\}
\mbox{ such that }\bigl |\sigma^{-1}(a)\bigr| = \mu_a, a=1,\ldots,K
\bigr\}.
\end{eqnarray*}

We are now equipped to introduce sets ${\mathcal{F}_{\mu\nu}^A}$ and
${\mathcal{F}_{\mu\nu}^\omega}$, which
describe all possible co-clusterings that can be induced from $A$ and
$\omega$ with respect to $(\mu, \nu) \in\Omega_m \times\Omega
_n$, and
to define the related notion of a support function.

\begin{definition}[(Sets ${\mathcal{F}_{\mu\nu}^A}$ and
${\mathcal{F}_{\mu\nu}^\omega}$ of admissible
co-clusterings)]\label{defn:clusterSets}
For fixed discrete probability distributions $\mu$ and $\nu$ over
$1,\ldots,K$, we define the sets ${\mathcal{F}_{\mu\nu}^A},
{\mathcal{F}_{\mu\nu}^\omega}\subset\mathbb{R}^{K
\times K}$ of all co-clustering matrices $A/ST$ and $\omega/\sigma
\tau
$, induced respectively by $(S,T) \in\mathcal{Q}_\mu^m \times
\mathcal
{Q}_\nu^n$ and $(\sigma, \tau) \in\mathcal{Q}_\mu\times\mathcal
{Q}_\nu$, as follows:
\begin{eqnarray*}
{\mathcal{F}_{\mu\nu}^A}&=& \bigl\{A/ST \in[0,1]^{K \times K}\dvtx S \in \mathcal{Q}_\mu^m, T \in\mathcal{Q}_\nu^n
\bigr\},
\\
{\mathcal{F}_{\mu\nu}^\omega}&=& \bigl\{\omega/\sigma\tau\in
[0,1]^{K \times K}\dvtx \sigma\in\mathcal {Q}_\mu, \tau\in
\mathcal{Q}_\nu \bigr\}.
\end{eqnarray*}
\end{definition}

\begin{definition}[(Support functions of ${\mathcal{F}_{\mu\nu
}^A}$ and ${\mathcal{F}_{\mu\nu}^\omega}$)]\label
{def:support}
Let $\mathcal{F} \subset\mathbb{R}^{K \times K}$ be nonempty and with
$\langle F, F' \rangle = \operatorname{tr}(F^T F')$. Its support function
$h_\mathcal{F}\dvtx \mathbb{R}^{K \times K} \rightarrow\mathbb{R} \cup
\{
+\infty\}$ is defined as $h_\mathcal{F}(\Gamma) = \sup_{F \in
\mathcal
{F}} \langle\Gamma, F \rangle$ for any $\Gamma\in\mathbb{R}^{K
\times K}$, whence
\begin{subequations}
\label{eq:supportFcns}
%
\begin{eqnarray}
h_{{\mathcal{F}_{\mu\nu}^A}}(\Gamma) & =& \max_{(S,T) \in\mathcal
{Q}_\mu^m \times\mathcal{Q}_\nu^n} \langle\Gamma,A/ST
\rangle,
\\
h_{{\mathcal{F}_{\mu\nu}^\omega}}(\Gamma) & =& \sup_{(\sigma,\tau
) \in\mathcal{Q}_\mu\times\mathcal{Q}_\nu} \langle\Gamma,\omega
/ \sigma\tau \rangle.
\end{eqnarray}
\end{subequations}
\end{definition}

We will show below that $\sup_{\Gamma\in[-1,1]^{K \times K}}
|h_{\smash{{\mathcal{F}_{\mu\nu}^A}}}(\Gamma) - h_{\smash
{{\mathcal{F}_{\mu\nu}^\omega}}}(\Gamma) |$ converges
in probability to zero at a rate of at least $n^{\smash{-1/4}}$, and
this result in turn gives rise to Theorem~\ref{th:learning}. To see
why, observe that for any $(\mu,\nu,\theta) \in\Phi$, the least
squares objective function of \eqref{eq:MSfunctional} can be expressed
using $h_{{\mathcal{F}_{\mu\nu}^A}}(\theta)$ as follows:
\begin{eqnarray*}
&&\min_{(S,T) \in\mathcal{Q}_\mu^m \times\mathcal{Q}_\nu^n}  \Biggl\{ \frac{1}{mn} \sum
_{i=1}^m \sum_{j=1}^n
\bigl( \theta_{S(i) T(j)}^2 - 2 \theta_{S(i)
T(j)}
A_{ij} + A_{ij}^2 \bigr) \Biggr\}
\\
&&\qquad = \sum_{a=1}^K \sum
_{b=1}^K \mu_a \nu_b
\theta_{ab}^2 - 2 \max_{(S,T) \in\mathcal{Q}_\mu^m \times\mathcal
{Q}_\nu^n} \langle
\theta, A/ST \rangle+ \frac{1}{mn} \sum_{i=1}^m
\sum_{j=1}^n A_{ij}
\\
&&\qquad = \sum_{a=1}^K \sum
_{b=1}^K \mu_a \nu_b
\theta_{ab}^2 - 2 h_{{\mathcal{F}_{\mu\nu}^A}
}(\theta) +
\frac{1}{mn} \sum_{i=1}^m \sum
_{j=1}^n A_{ij}.
\end{eqnarray*}
As we prove in Appendix \ref{sec:proof.th2}, this line of argument
establishes the following.

\begin{lemma}\label{lem:likelihood_support_fcns}
For any $(\mu,\nu,\theta) \in\Phi$, the difference between the least
squares objective function of \eqref{eq:MSfunctional} and the $L^2$
risk $R_\omega$ is equal to
\[
2 \bigl( h_{{\mathcal{F}_{\mu\nu}^\omega}}(\theta) - h_{{\mathcal
{F}_{\mu\nu}^A}}(\theta) \bigr) +
\frac
{1}{mn}\sum_{i=1}^m \sum
_{j=1}^n A_{ij}^2 -
\int_{[0,1]^2} \omega(x,y)^2 \,dx \,dy,
\]
and the difference between the profile likelihood function of \eqref
{eq:PLfunctional} and $L_\omega$ is $B(\theta) ( h_{{\mathcal
{F}_{\mu\nu}^A}}( \Gamma
_\theta) - h_{{\mathcal{F}_{\mu\nu}^\omega}}(\Gamma_\theta) )$
whenever $0 < \theta_{ab}
<1$ for all $a,b=1,\ldots,K$, with $\Gamma_\theta\in[-1,1]^{\smash{K
\times K}}$ given element-wise by $(\Gamma_\theta)_{ab} = \log
(\theta
_{ab} / (1-\theta_{ab})) / B(\theta)$.
\end{lemma}

\subsection{A general result on consistency of co-clustering}

From Lemma~\ref{lem:likelihood_support_fcns} we see that closeness of
$h_{\smash{{\mathcal{F}_{\mu\nu}^A}}}$ to $h_{\smash{{\mathcal
{F}_{\mu\nu}^\omega}}}$ implies closeness (up to
constant terms) of the least squares objective function of \eqref
{eq:MSfunctional} to the $L^2$ risk $R_{\smash{\omega}}(\phi)$, and of
the profile likelihood of \eqref{eq:PLfunctional} to the average
Kullback--Leibler divergence of $\omega_\phi(x,y)$ from the generative
$\omega(x,y)$. Equipped with this motivation, we now state our main
technical result, which serves to establish the rate of convergence
$\mathcal{O}_P (n^{\smash{-1/4}} )$ in Theorem~\ref{th:learning}.
Its proof follows in Section~\ref{subsec:sketch} below.

\begin{theorem}\label{th:main}
Let $A \in\{0,1\}^{m \times n}$ be a separately exchangeable array
generated by some $\omega$ in accordance with Definition~\ref
{model:general}. Then for each $K$ and each ratio $m/n$, there exists a
universal constant $C$ such that as $n \rightarrow\infty$,
\[
\mathbb{P} \biggl( \max_{(\mu,\nu) \in\Omega_m \times\Omega_n} \Bigl\{ \sup
_{\Gamma\in[-1,1]^{K \times K}}\bigl | h_{{\mathcal{F}_{\mu\nu
}^A}}(\Gamma) - h_{{\mathcal{F}_{\mu\nu}^\omega}
}(\Gamma)
\bigr| \Bigr\} \geq\frac{C}{n^{1/4}} \biggr) = o(1).
\]
\end{theorem}

The support functions $h_{\smash{{\mathcal{F}_{\mu\nu}^A}}}$ and
$h_{\smash{{\mathcal{F}_{\mu\nu}^\omega}}}$ also
have a geometric interpretation: for any fixed $\Gamma\in\mathbb
{R}^{K \times K}$, they define the supporting hyperplanes of the sets
${\mathcal{F}_{\mu\nu}^A}$ and ${\mathcal{F}_{\mu\nu}^\omega}$
in the direction specified by $\Gamma$. Each
supporting hyperplane is induced by a point in ${\mathcal{F}_{\mu\nu
}^A}$, or in the
closure of ${\mathcal{F}_{\mu\nu}^\omega}$ respectively; these
points are extremal in that they
cannot be written as a convex combination of any other points in their
respective sets. Evidently, it is only the extreme points which
determine convergence properties for the risk functionals considered
here. Equivalently, for any fixed parameter triple $\phi\in\Phi$, the
values of these functionals depend only on the maximizing choices of
$(S,T)$ or $(\sigma,\tau)$.

Formally, Theorem~\ref{th:main} has the following geometric interpretation:
%
\begin{corollary}\label{corr:convHulls}
The result of Theorem \ref{th:main}
is equivalent to the following: The
Hausdorff distance between the convex hulls of ${\mathcal{F}_{\mu\nu
}^A}$ and ${\mathcal{F}_{\mu\nu}^\omega}$ is
$\mathcal{O}_P ( n^{-1/4} )$.
\end{corollary}

\begin{pf}
Consider $\mathcal{F}, \mathcal{F}' \subset\mathbb{R}^{K \times K}$,
and denote by $\|F\| = \sqrt{\operatorname{tr}(F^TF)}$ the Frobenius
norm (i.e., the Hilbert--Schmidt metric on $\mathbb{R}^{K \times K}$
induced by $\langle\cdot,\cdot \rangle$). The Hausdorff distance between
$\mathcal{F}$ and $\mathcal{F}'$, based on the metric $\|\cdot\|$,
is then
\[
d_{\mathrm{Haus}} \bigl( \mathcal{F},\mathcal{F}' \bigr) = \max
\Bigl\{ \sup_{F \in\mathcal{F}} \Bigl\{ \inf_{{F'} \in
{\mathcal{F}'}}
\bigl\|F-F'\bigr\| \Bigr\}, \sup_{F' \in\mathcal{F}'} \Bigl\{ \inf
_{F \in\mathcal{F}} \bigl\|F-F'\bigr\| \Bigr\} \Bigr\}.
\]
This measures the maximal shortest distance between any two elements of
$\mathcal{F}$ and~$\mathcal{F}'$. If these subsets of $\mathbb{R}^{K
\times K}$ are furthermore nonempty and bounded, then the Hausdorff
distance between their convex hulls $\operatorname{conv} (\mathcal
{F} )$ and $\operatorname{conv} (\mathcal{F}' )$ can be
expressed in terms of their support functions $h_{\mathcal{F}},
h_{\mathcal{F}'}$,
\[
d_{\mathrm{Haus}} \bigl( \operatorname{conv} (\mathcal{F} ), \operatorname{conv}
\bigl(\mathcal{F}' \bigr) \bigr) = \sup_{\Gamma\in\mathbb
{R}^{K\times K} \dvtx \|\Gamma\| = 1} \bigl|
h_{\mathcal{F}}(\Gamma) - h_{\mathcal{F}'}(\Gamma) \bigr|;
\]
see, for example, \citet{schneider1993convex}, as applied to the convex
hulls of the closures of $\mathcal{F}$ and of $\mathcal{F}'$. In this
way, $d_{\mathrm{Haus}}(\cdot,\cdot)$ is a natural measure of distance
between two convex bodies. Recalling the equivalence of norms on
$\mathbb{R}^{K^2}$, we see that
\begin{eqnarray*}
\sup_{\|\Gamma\| = 1} \bigl| h_{\mathcal{F}}(\Gamma) - h_{\mathcal
{F}'}(
\Gamma) \bigr| &\leq&\sup_{\Gamma\in[-1,1]^{K \times K}} \bigl| h_{\mathcal
{F}}(\Gamma) -
h_{\mathcal{F}'}(\Gamma) \bigr|\\
& \leq& K \sup_{\|\Gamma\| = 1} \bigl|
h_{\mathcal{F}}(\Gamma) - h_{\mathcal
{F}'}(\Gamma) \bigr|.
\end{eqnarray*}
Since Theorem~\ref{th:main} holds for $\sup_{\Gamma\in[-1,1]^{K
\times K}} | h_{\smash{{\mathcal{F}_{\mu\nu}^A}}}(\Gamma) -
h_{\smash{{\mathcal{F}_{\mu\nu}^\omega}}}(\Gamma) |$,
the leftmost inequality implies that it also holds for $\sup_{\|\Gamma
\| = 1} | h_{\smash{{\mathcal{F}_{\mu\nu}^A}}}(\Gamma) - h_{\smash
{{\mathcal{F}_{\mu\nu}^\omega}}}(\Gamma) |$. Now
suppose instead that Theorem~\ref{th:main} holds for $\sup_{\|\Gamma
\|
= 1} | h_{\smash{{\mathcal{F}_{\mu\nu}^A}}}(\Gamma) - h_{\smash
{{\mathcal{F}_{\mu\nu}^\omega}}}(\Gamma) |$; by the
rightmost inequality, it then also holds for $ K^{-1} \sup_{\Gamma\in
[-1,1]^{K \times K}} | h_{\smash{{\mathcal{F}_{\mu\nu}^A}}}(\Gamma
) - h_{\smash{{\mathcal{F}_{\mu\nu}^\omega}
}}(\Gamma) |$. Thus the result of Theorem~\ref{th:main} is equivalent
to the statement that
\[
\max_{(\mu,\nu) \in\Omega_m \times\Omega_n} d_{\mathrm{Haus}} \bigl( \operatorname{conv}
\bigl( {\mathcal{F}_{\mu\nu}^A}\bigr),\operatorname{conv}
\bigl( {\mathcal{F}_{\mu
\nu}^\omega}\bigr) \bigr) =
\mathcal{O}_P \bigl( n^{-1/4} \bigr).
\]
\upqed\end{pf}

This geometric interpretation is helpful in relating our work to a
series of papers by \citeauthor{borgs2006graph} (\citeyear{borgs2006graph,borgs2008convergent,borgs2007convergent}), which explore dense graph limits in depth and
statistical applications thereof. Very broadly speaking,
\citet{borgs2008convergent}, Theorem~2.9 and
\citet{borgs2007convergent}, Theorem~4.6,
analyze sets termed quotients, which resemble $\bigcup_{\mu,\nu}
{\mathcal{F}_{\mu\nu}^A}$
and $\bigcup_{\mu,\nu} {\mathcal{F}_{\mu\nu}^\omega}$. The authors
show convergence of these sets
in the Hausdorff metric at rate $\mathcal{O} ( \log^{\smash{-1/2}} n
)$, based on a distance termed the cut metric, and detail
implications that can also be related to those of \citet{bickelmethod}.

In fixing $\mu$ and $\nu$ through our $M$-estimators, we are studying
what Borgs et al. term the microcanonical quotients. Because our
results require only convergence of the closed convex hulls of
${\mathcal{F}_{\mu\nu}^A}$
and ${\mathcal{F}_{\mu\nu}^\omega}$, we are able to obtain an
exponentially faster bound on the
rate of convergence.

\subsection{Interpreting convergence of blockmodel estimates}

Recall that the $M$-estimators of Theorem~\ref{th:learning} each involve
an optimization over the set ${\mathcal{F}_{\mu\nu}^A}$ by way of
its support function,
which in turn represents its convex hull. Suppose that $\smash{\hat
{\phi
}} \equiv(\hat{\mu},\hat{\nu}, \smash{\hat{\theta}})$ optimizes either
objective function in Theorem~\ref{th:learning}. Then the following
corollary of Theorem~\ref{th:learning} shows that $\hat{\phi}$ is
interpretable, in that there will exist a partition $\hat{\sigma},
\hat
{\tau}$ of $\omega$ yielding co-clusters of equal size and
asymptotically equivalent connectivity.

\begin{corollary}\label{le:interpretable}
Let $\hat{\phi} = (\hat{\mu},\hat{\nu}, \hat{\theta})$ minimize the
least squares criterion of~\eqref{eq:MSfunctional}. Then there exists
some pair $(\hat{\sigma}, \hat{\tau}) \in\mathcal{Q}_{\hat{\mu}}
\times\mathcal{Q}_{\hat{\nu}}$ such that
\[
\sum_{a=1}^K \sum
_{b=1}^K \hat{\mu}_a \hat{
\nu}_b \biggl| \frac{(\omega
/\hat{\sigma}\hat{\tau})_{ab}}{\hat{\mu}_a \hat{\nu}_b} - \hat {\theta}_{ab}
\biggr|^2 = \mathcal{O}_P \bigl(n^{-1/4} \bigr).
\]
Similarly, if $\hat{\phi} = (\hat{\mu},\hat{\nu}, \hat{\theta})$
maximizes the profile likelihood criterion of \eqref{eq:PLfunctional}
and $\phi^* = \operatorname{argmax}_{\phi\in\Phi} L_\omega(\phi)$
exists, then there is some $(\hat{\sigma}, \hat{\tau}) \in\mathcal
{Q}_{\hat{\mu}} \times\mathcal{Q}_{\hat{\nu}}$ with
\[
\frac{1}{B(\phi^*) + B(\hat{\phi})} \sum_{a=1}^K \sum
_{b=1}^K \hat{\mu}_a \hat{
\nu}_b D \biggl(\frac{(\omega/\hat{\sigma}\hat{\tau
})_{ab}}{\hat{\mu}_a \hat{\nu}_b} \Big\Vert\hat{\theta
}_{ab} \biggr) = \mathcal{O}_P \bigl(n^{-1/4}
\bigr),
\]
where $D(p \Vert p') = p \log(p/p') + (1-p) \log[(1-p)/(1-p')]
\geq0$ is the Kullback--Leibler divergence of a $\operatorname
{Bernoulli}(p')$ distribution from a $\operatorname{Bernoulli}(p)$ one.
\end{corollary}

\begin{pf}
We show the latter result; parallel arguments yield the former. Since
$\omega_{\hat{\phi}}(x,y) = \smash{\hat{\theta}}_{\smash
{F^{-1}_{\hat
{\mu}}(x) F^{-1}_{\hat{\nu}}(y)}}$ for the co-blockmodel, by letting
$\sigma$ and $\tau$ satisfy $\sigma(x) = F^{\smash{-1}}_{\hat{\mu
}}(\pi
_1(x))$ and $\tau(y) = F^{\smash{-1}}_{\hat{\nu}}(\pi_2(y))$ we may
express $L_\omega(\smash{\hat{\phi}})$ as
\begin{eqnarray*}
&&\sup_{(\sigma,\tau) \in\mathcal{Q}_{\hat{\mu}} \times\mathcal
{Q}_{\hat
{\nu}}} \sum_{a=1}^K
\sum_{b=1}^K \int_{\sigma^{-1}(a)\times\tau
^{-1}(b)}
\bigl\{ \omega(x,y) \log\hat{\theta}_{ab} \\
&&\hspace*{153pt}{}+ \bigl[ 1-\omega(x,y)
\bigr] \log(1-\hat{\theta}_{ab}) \bigr\} \,dx \,dy.
\end{eqnarray*}
Thus, for any $\varepsilon>
 0$, there exists some choice of $(\hat
{\sigma
}, \hat{\tau}) \in\mathcal{Q}_{\smash{\hat{\mu}}} \times
\mathcal
{Q}_{\smash{\hat{\nu}}}$ such that
\[
L_\omega(\hat{\phi}) - \varepsilon\leq\sum_{a=1}^K
\sum_{b=1}^K \bigl\{ (\omega/\hat{\sigma}
\hat{\tau})_{ab} \log\hat{\theta}_{ab} + \bigl[ \hat{
\mu}_a \hat{\nu}_b - (\omega/\hat{\sigma} \hat{
\tau})_{ab} \bigr] \log(1-\hat{\theta}_{ab}) \bigr\}.
\]
If we now take $\hat{\theta}^{(\omega)}_{\smash{ab}} = (\omega
/\hat
{\sigma}\hat{\tau})_{ab} / (\hat{\mu}_a \hat{\nu}_b)$ for
$a,b=1,\ldots
,K$, we see by a similar argument that since $L_\omega(\phi^*) = \max_{\phi\in\Phi} L_\omega(\phi)$, we have in turn that
\begin{eqnarray*}
L_\omega \bigl(\phi^* \bigr) & \geq& L_\omega \bigl( \bigl(\hat{
\mu}, \hat{\nu}, \hat{\theta}^{(\omega)} \bigr) \bigr)
\\
& \geq&\sum_{a=1}^K \sum
_{b=1}^K \bigl\{ (\omega/\hat{\sigma} \hat{\tau
})_{ab} \log\hat{\theta}^{(\omega)}_{ab} + \bigl[\hat{
\mu}_a \hat{\nu}_b - (\omega/\hat{\sigma} \hat{
\tau})_{ab} \bigr] \log \bigl(1- \hat{\theta}^{(\omega)}_{ab}
\bigr) \bigr\}.
\end{eqnarray*}
Expanding $D (\hat{\theta}^{(\omega)}_{ab} \Vert
\hat{\theta}_{ab} )$ in accordance with its definition, we then see that
\[
0 \leq\sum_{a=1}^K \sum
_{b=1}^K \hat{\mu}_a \hat{
\nu}_b D \bigl(\hat{\theta}^{(\omega)}_{ab} \Vert\hat{\theta}_{ab} \bigr) \leq L_\omega \bigl(\phi^*
\bigr)- L_\omega(\smash{\hat{\phi}}) + \varepsilon.
\]
Choosing $\varepsilon= o(n^{\smash{-1/4}})$ and applying Theorem~\ref
{th:learning} completes the proof.
\end{pf}

Corollary~\ref{le:interpretable} ensures that co-blockmodel fits remain
interpretable, even in the setting of model misspecification. It
establishes that the identification of co-clusters in an observed
exchangeable binary array $A$ indicates with high probability the
existence of co-clusters of equal size and asymptotically equivalent
connectivity in the underlying generative process $\omega$.

\section{\texorpdfstring{Proof of Theorem \protect\ref{th:main}}
{Proof of Theorem 4.1}}
\label{subsec:sketch}

Our proof strategy is inspired by \citet{borgs2008convergent} and
adapts certain of its tools, but also requires new techniques in order
to attain polynomial rates of convergence. Most significantly, we do
not use the Szemer\'{e}di regularity lemma, which typically features
strongly in the graph-theoretic literature, and provides a means of
partitioning any large dense graph into a small number of regular
clusters. Results in this direction are possible, but instead we use a
Rademacher complexity bound for $U$-statistics adapted from \citet
{clemenccon2008ranking}, allowing us to achieve the improved rates of
convergence described above.

\subsection{Establishing pointwise convergence}
\label{sec:pwConv}

The main step in proving Theorem~\ref{th:main} is to establish
pointwise convergence of $h_{{\mathcal{F}_{\mu\nu}^A}}(\Gamma)$ to
$h_{{\mathcal{F}_{\mu\nu}^\omega}}(\Gamma)$ for
any fixed~$\Gamma$. We do this through Proposition~\ref
{prop:hyperplane} below, after which we may apply it to a union bound
over a covering of all $\Gamma\in[-1,1]^{K \times K}$ to deduce the
result of Theorem~\ref{th:main}. Appendix \ref{sec:proof.th3} provides
a formal statement and proof of this argument, along with proofs of all
supporting lemmas.

\begin{proposition}[{[Pointwise convergence of $h_{{\mathcal{F}_{\mu\nu
}^A}}(\Gamma)$ to
$h_{{\mathcal{F}_{\mu\nu}^\omega}}(\Gamma)$]}]\label{prop:hyperplane}
Assume the setting of Theorem~\ref{th:main}, fixing $m = \rho n$. Then
there exist constants $C_K,n_K$ such that, given any $\Gamma\in
[-1,1]^{K \times K}$, $\mu$, $\nu$, $\omega$, and $A \in\{0,1\}^{m
\times n}$ generated from $\omega$, it holds for all $n \geq n_K$ that
\[
\mathbb{P} \biggl( \bigl| h_{{\mathcal{F}_{\mu\nu}^A}}(\Gamma) -
h_{{\mathcal{F}_{\mu\nu}^\omega}}(\Gamma) \bigr|
\geq\frac{C_K }{n^{1/4}} \biggr) \leq2e^{-\sqrt{n} [2\rho/(\rho
+1)]} \bigl[1+o(1) \bigr].
\]
\end{proposition}

\begin{pf}
To obtain the claimed result, we must establish lower and upper bounds
on the support function $h_{{\mathcal{F}_{\mu\nu}^A}}(\Gamma)$ that
show its convergence
to $h_{{\mathcal{F}_{\mu\nu}^\omega}}(\Gamma)$ at rate $\mathcal
{O}_P (n^{\smash
{-1/4}} )$. Recalling the definitions of $h_{{\mathcal{F}_{\mu\nu
}^A}}(\Gamma)$ and
$h_{{\mathcal{F}_{\mu\nu}^\omega}}(\Gamma)$ in \eqref
{eq:supportFcns}, we first require a
statement of Lipschitz conditions on $\langle\Gamma,A/ST \rangle$
and $\langle\Gamma,\omega/\sigma\tau \rangle$. Its proof follows
by direct inspection.

\begin{lemma} \label{le:lipschitz}
Define for measurable mappings $\sigma, \sigma'$ over $[0,1]$ the metric
\[
d_\mathrm{Ham} \bigl(\sigma,\sigma' \bigr) = \int
_{[0,1]} 1 \bigl\{\sigma(x) \neq\sigma'(x) \bigr\}
\,dx,
\]
and analogously the standard Hamming distance for sequences, with
respect to normalized counting measure. Then for any $\Gamma\in
[-1,1]^{K \times K}$ and $A,A' \in[0,1]^{m \times n}$, with
$(S,T,\omega,\sigma,\tau)$ as defined in Section~\ref{sec:notation}, we
have that:
\begin{longlist}[(1)]
\item[(1)] $|\langle\Gamma,A/ST \rangle - \langle\Gamma,A/S'T'
\rangle| \leq2 [
d_{\mathrm{Ham}}(S,S') / m + d_{\mathrm{Ham}}(T,T') / n ]$;
\item[(2)] $|\langle\Gamma,\omega/\sigma\tau \rangle - \langle
\Gamma, \omega/\sigma '\tau' \rangle| \leq2 [ d_{\mathrm
{Ham}}(\sigma,\sigma') + d_{\mathrm
{Ham}}(\tau,\tau') ]$;
\item[(3)] $|\langle\Gamma,A/ST \rangle - \langle\Gamma,A'/ST
\rangle| \leq1/(mn)$ if
$A,A'$ differ by a single entry.
\end{longlist}
\end{lemma}

In conjunction with McDiarmid's inequality, these Lipschitz conditions
yield the following lower bound on $h_{{\mathcal{F}_{\mu\nu
}^A}}(\Gamma)$, proved in
Appendix \ref{le:lowerpf}.

\begin{lemma}[{[Lower bound on $h_{{\mathcal{F}_{\mu\nu}^A}}(\Gamma
)$]}]\label{le:lower}
Assume the setting of Theorem~\ref{th:main}. Then there exist constants
$C{}^{\prime}_K, n{}^{\prime}_K $ such that, given any $\Gamma\in
[-1,1]^{K \times K}, \mu, \nu, \omega$, and $A \in\{0,1\}^{\rho n
\times n}$ generated from $\omega$, for all $n \geq n_K '$,
\[
\mathbb{P} \biggl(h_{{\mathcal{F}_{\mu\nu}^\omega}}(\Gamma) - h_{{\mathcal{F}_{\mu\nu}^A}}(\Gamma) \geq
\frac{C_K
'}{n^{1/4}} \biggr) \leq2e^{-\sqrt{n} [2\rho/(\rho+1)]} \bigl[1+o(1) \bigr].
\]
\end{lemma}

The upper bound comes by way of Rademacher complexity arguments. The
remainder of this section and Appendix \ref{sec:proof.th3} is devoted
to its proof.

\begin{lemma}[{[Upper bound on $h_{{\mathcal{F}_{\mu\nu}^A}}(\Gamma
)$]}]\label{le:upper}
Assume the setting of Theorem~\ref{th:main}. Then there exist constants
$C^{\prime\prime}_K, n^{\prime\prime}_K $ such that, given any
$\Gamma\in\break
[-1,1]^{K \times K}, \mu, \nu, \omega$ and $A \in\{0,1\}^{\rho n
\times n}$ generated from $\omega$, for all $n \geq n_K ''$,
\[
\mathbb{P} \biggl(h_{{\mathcal{F}_{\mu\nu}^A}}(\Gamma) - h_{{\mathcal{F}_{\mu\nu}^\omega}}(\Gamma) \geq
\frac{C_K
''}{n^{1/4}} \biggr) \leq2e^{-\sqrt{n} [2\rho/(\rho+1)]} \bigl[1+o(1) \bigr].
\]
\end{lemma}

Proposition~\ref{prop:hyperplane} now follows simply by combining
Lemmas \ref{le:lower} and \ref{le:upper}.
\end{pf}

\subsection{\texorpdfstring{Establishing an upper bound on $h_{{\mathcal{F}_{\mu\nu}^A}}(\Gamma)$}
{Establishing an upper bound on h F mu nu A (Gamma)}}
\label{sec:K2Case}

Lemma~\ref{le:upper} represents the main technical hurdle in obtaining
the polynomial rate of convergence given in Theorems \ref{th:learning}
and \ref{th:main}. To illustrate the main ideas as clearly as possible,
we will introduce our Rademacher complexity arguments below for the
case $K = 2$, deferring the necessary generalizations to Appendix \ref
{sec:proof.th3}.

We first define $W \in[0,1]^{m\times n}$ with reference to
Definition~\ref{model:general} as
\[
W_{ij} = \omega(\xi_i,\zeta_j),\qquad i \in1,
\ldots,m, j \in1,\ldots,n;
\]
and then define, in direct analogy to $h_{{\mathcal{F}_{\mu\nu
}^A}}(\Gamma)$,
\[
h_{\mathcal{F}_{\mu\nu}^W}(\Gamma) = \max_{(S,T) \in\mathcal
{Q}_\mu^m \times\mathcal{Q}_\nu^n} \langle\Gamma,W/ST
\rangle = \max_{(S,T)
\in\mathcal{Q}_\mu^m \times\mathcal{Q}_\nu^n} \Biggl\{ \frac
{1}{mn} \sum
_{i=1}^m \sum_{j=1}^n
W_{ij} \Gamma_{S(i)T(j)} \Biggr\}.
\]

The matrix $W$ serves as an empirical realization of the mapping
$\omega
$, with its support function $h_{\mathcal{F}_{\mu\nu}^W}(\Gamma)$
defined with respect to co-blockmodel partitions $(S,T) \in\mathcal
{Q}_\mu^m \times\mathcal{Q}_\nu^n$. As proved in Appendix \ref
{le:basicpf}, Lemma~\ref{le:basic} enables us to bound $ |
h_{{\mathcal{F}_{\mu\nu}^A}
}(\Gamma) - \mathbb{E} h_{\mathcal{F}_{\mu\nu
}^W}(\Gamma
) |$ using the Lipschitz conditions in Lemma~\ref{le:lipschitz}.

\begin{lemma}\label{le:basic}
Fix some measurable $\omega\dvtx [0,1]^2\rightarrow[0,1]$, with $W \in
[0,1]^{m \times n}$ generated by $\omega$ and $A \in\{0,1\}^{m \times
n}$ generated by $W$, and some $\Gamma\in[-1,1]^{K \times K}$. Then
for any $\varepsilon> 0$,
%
\begin{equation}
\label{eq:basic3} \mathbb{P} \bigl( \bigl| h_{{\mathcal{F}_{\mu\nu
}^A}}(\Gamma) -
\mathbb{E} h_{\mathcal{F}_{\mu\nu}^W}(\Gamma)\bigr | \geq2\varepsilon \bigr) \leq2e^{-2mn\varepsilon^2/(m+n)}
+ 2K^{m+n} e^{-2mn\varepsilon^2}.
\end{equation}
\end{lemma}

Having bounded $ | h_{\smash{{\mathcal{F}_{\mu\nu}^A}}}(\Gamma) -
\mathbb
{E} h_{\mathcal{F}_{\smash{\mu\nu}}^W}(\Gamma) |$,\vspace*{1pt} we must
upper-bound $\mathbb{E} h_{\mathcal{F}_{\smash{\mu
\nu
}}^W}(\Gamma)$ in terms of $ h_{\smash{{\mathcal{F}_{\mu\nu
}^\omega}}}(\Gamma)$. We do this in
a series of steps, first bounding $\mathbb{E}
h_{\mathcal
{F}_{\smash{\mu\nu}}^W}(\Gamma)$ using a result adapted from \citet
{alon2003random} and proved in Appendix \ref{le:alonpf}.

\begin{lemma}\label{le:alon}
Let $\mathcal{I}$ and $\mathcal{J}$ be sets of deterministic size,
whose elements are sampled without replacement from $1,\ldots,m$ and
$1,\ldots,n$. Let $W$ be generated as in Lemma~\ref{le:basic}, and fix
$\Gamma\in[-1,1]^{K \times K}$. Given $W,\mathcal{I},\mathcal{J}$ and
$(Q, R) \in\mathcal{Q}_{\mu}^m \times\mathcal{Q}_{\nu}^n$, let
$\hat
{S}^R \equiv\hat{S}^{R,\mathcal{J},W}$ and $\hat{T}^Q \equiv\hat
{T}^{Q,\mathcal{I},W}$ denote partitions satisfying
%
\begin{eqnarray}
\label{eq:hatS} \hat{S}^R & =& \mathop{\operatorname{argmax}}_{S \in\mathcal
{Q}_\mu^m}
\Biggl\{ \sum_{i=1}^m \sum
_{j \in\mathcal{J}} W_{ij} \Gamma_{S(i)R(j)} \Biggr\},
\\
\label{eq:hatT} \hat{T}^Q & = &\mathop{\operatorname{argmax}}_{T \in\mathcal
{Q}_\nu^n}
\Biggl\{ \sum_{i \in
\mathcal{I}} \sum
_{j=1}^n W_{ij} \Gamma_{Q(i)T(j)}
\Biggr\}.
\end{eqnarray}
Then
%
\begin{eqnarray}
\label{eq:alon} \mathbb{E} h_{\mathcal{F}_{\mu\nu
}^W}(\Gamma) &\leq&
\mathbb{E} \Bigl(\max_{(Q,R) \in\mathcal{Q}_\mu
^m\times
\mathcal{Q}_{\nu}^n} \bigl\langle
\Gamma,W/ \hat{S}^R\hat{T}^Q \bigr\rangle \Bigr)
\nonumber
\\[-8pt]
\\[-8pt]
\nonumber
&&{}+ K
\sqrt{2\pi} \bigl( | \mathcal{I}|^{-{1}/{2}} + |\mathcal{J}|^{-
{1}/{2}}
\bigr).
\end{eqnarray}
\end{lemma}

To bound the right-hand side of \eqref{eq:alon} relative to
$h_{{\mathcal{F}_{\mu\nu}^\omega}
}(\Gamma)$, we will introduce an additional construction comprising
several steps. Specifically, for fixed $(Q,R)$ and $\Gamma$, we will
define function classes $\mathcal{Q}_U$ and $\mathcal{Q}_V$, and a
random functional $G_{\sigma\tau}$ which approximates $\langle\Gamma
,W/\hat {S}^R\hat{T}^Q \rangle$ for some $(\hat{\sigma},\hat{\tau
}) \in\mathcal{Q}_U
\times\mathcal{Q}_V$. By a Rademacher complexity argument, $G_{\hat
{\sigma}\hat{\tau}}$ will concentrate for all $(Q,R)$ near its
expectation, which itself will be bounded by $h_{\mathcal{F}_{\mu\nu
}^\omega}(\Gamma)$.

For the case $K=2$, define $U$ by
\[
U(x) = \sum_{j \in\mathcal{J}} \omega(x,\zeta_j) (
\Gamma_{1R(j)} - \Gamma_{2R(j)}).
\]
It follows that
\[
\hat{S}^R = \mathop{\operatorname{argmax}}_{S \in\mathcal{Q}_\mu^m} \sum
_{i=1}^m U(\xi_i) 1 \bigl\{S(i)=1
\bigr\},
\]
and so $\hat{S}^R$ will assign to class $1$ the $\mu_1m$ largest
elements of $U(\xi_1),\ldots,U(\xi_m)$. If $U$ is invertible, this set
can be written $\{\xi_i\dvtx U(\xi_i) < t\}$ for some $t$. To treat
noninvertible $U$, define $\mathcal{Q}_U$ to be the class of functions
$\{1_u\dvtx u\in[0,1]\}$, with $1_u$ a one-sided interval on the range of
$U$ with lexicographic ``tie-breaking'':
\[
1_u(x) = \cases{ 2, &\quad $\mbox{if either $U(x) < U(u)$, or $U(x)=U(u)$
and $x < u$;}$ \vspace*{2pt}
\cr
1, &\quad $ \mbox{if either $U(x) > U(u)$, or
$U(x)=U(u)$ and $x \geq u$.}$}
\]
Then there exists $\hat{\sigma} \in\mathcal{Q}_U$ such that $\hat
{S}^R$ can be chosen to satisfy
\[
\hat{S}^R(i) = \hat{\sigma}(\xi_i),\qquad i=1,\ldots,m.
\]

Let $V$ denote a function defined analogously to $U$ as follows:
\[
V(y) = \sum_{i \in\mathcal{I}} \omega(\xi_i,y) (
\Gamma_{Q(i)1} - \Gamma_{Q(i)2}),
\]
and likewise define $\mathcal{Q}_V$ so that there exists $\hat{\tau}
\in\mathcal{Q}_V$ such that $\hat{T}^Q$ can be chosen to satisfy
\[
\hat{T}^Q(j) = \hat{\tau}(\zeta_j),\qquad j=1,\ldots,n.
\]

We are now ready to define $G_{\sigma\tau}$. Given any $\sigma\in
\mathcal{Q}_U$ and $\tau\in\mathcal{Q}_V$, let
\[
G_{\sigma\tau}(\xi,\zeta)  = \frac{1}{mn} \sum
_{i \in\overline
{\mathcal{I}}} \sum_{j \in\overline{\mathcal{J}}} \omega(
\xi_i,\zeta_j)\Gamma_{\sigma(\xi_i)\tau(\zeta_j)},
\]
where $\overline{\mathcal{I}}$ is the complement of $\mathcal{I}$ in
$\{
1,\ldots,m\}$, and $\overline{\mathcal{J}}$ the complement of
$\mathcal
{J}$ in $\{1,\ldots,n \}$. Comparing $G_{\sigma\tau}$ to Lemma~\ref
{le:alon}, we see that $G_{\hat{\sigma}\hat{\tau}}$ well approximates
$\langle\Gamma,W/\hat{S}^R\hat{T}^Q \rangle$ whenever $|\mathcal
{I}|$ and
$|\mathcal{J}|$ are small; and indeed, we will later set $|\mathcal{I}|
= |\mathcal{J}| = n^{1/2}$ in order to obtain an upper bound for
$h_{{\mathcal{F}_{\mu\nu}^A}}(\Gamma) - h_{\mathcal{F}_{\mu\nu
}^\omega}(\Gamma)$.

By construction, the random classes $\mathcal{Q}_U$ and $\mathcal{Q}_V$
are independent of the random variables $\{\xi_i\}_{i \in\overline
{\mathcal{I}}}$ and $\{\zeta_j\}_{i \in\overline{\mathcal{J}}}$
appearing in the summand of $G_{\sigma\tau}$. As a result, we may bound
the deviation $\delta_{UV}$ of $G_{\sigma\tau}$ from its expectation,
\[
\delta_{UV} = \sup_{(\sigma,\tau) \in\mathcal{Q}_U\times\mathcal
{Q}_V} \bigl\llvert
G_{\sigma\tau}(\xi,\zeta) -\mathbb{E} \bigl(
G_{\sigma\tau}(\xi,\zeta) \vert U,V \bigr) \bigr\rrvert,
\]
using Rademacher complexity results for $U$-statistics due to \citet
{hoeffding1963probability} and \citet{clemenccon2008ranking}, Lemma A.1, applied to the class of one-sided
interval functions.

\begin{lemma} \label{le:rademacher1.easy}
Assume the setting of Lemma~\ref{le:alon}, and set $\ell= \min
(m-|\mathcal{I}|, n -|\mathcal{J}|)$. Then the deviation $\delta_{UV}$
of $G_{\sigma\tau}$ from its expectation satisfies
\[
\mathbb{E} \Bigl(\max_{(Q,R) \in\mathcal{Q}_\mu^m
\times\mathcal{Q}_\nu^n}
\delta_{UV} \Bigr) \leq4 \sqrt{\frac{(|\mathcal{I}| + |\mathcal
{J}|)\log K + 2{K\choose
2}\log(\ell+1) + \log2}{2 \ell}}.
\]
\end{lemma}

Lemma~\ref{le:rademacher1.easy} is proved in Appendix \ref
{le:rademacher1.easypf} to hold for arbitrary $K$, under the
appropriate generalization of $\mathcal{Q}_U, \mathcal{Q}_V$, and
quantities that depend on them.

Similarly, we may bound $\delta_U$, defined for $K = 2$ as the maximum
discrepancy between the expected and empirical class frequency in
$\mathcal{Q}_U$,
\[
\delta_U = \sup_{\sigma\in\mathcal{Q}_U} \Biggl\{ \max
_{1 \leq a \leq K}\Biggl | \bigl|\sigma^{-1}(a)\bigr| - \frac{1}{m}\sum
_{i=1}^m 1 \bigl\{\sigma(
\xi_i)=a \bigr\} \Biggr| \Biggr\},
\]
with $\delta_V$ defined mutatis mutandis. We then have the following
result, proved for arbitrary $K$ (with appropriate redefinitions of
$\delta_U, \delta_V$) in Appendix \ref{le:rademacher2.easypf}.

\begin{lemma}
\label{le:rademacher2.easy}
Assume the setting of Lemma~\ref{le:alon}. Then
\begin{eqnarray*}
\mathbb{E} \Bigl( \max_{R \in\mathcal{Q}_\nu^n} \delta
_U \Bigr) & \leq&4 \sqrt{\frac{(|\mathcal{J}|+1)\log K + {K\choose
2} \log(m+1)
+ \log2}{2m}},
\\
\mathbb{E} \Bigl( \max_{Q \in\mathcal{Q}_\mu^m}
\delta_V \Bigr) & \leq&4 \sqrt{\frac{(|\mathcal{I}|+1)\log K +
{K\choose2} \log(n+1)
+ \log2}{2n}}.
\end{eqnarray*}
\end{lemma}

We state and prove a final auxiliary lemma prior to the proof of
Lemma~\ref{le:upper}.

\begin{lemma}
\label{le:final}
Assume the setting of Lemma~\ref{le:alon}. Then
\begin{eqnarray*}
&&\mathbb{E} \Bigl( \max_{(Q,R) \in\mathcal{Q}_\mu^m
\times\mathcal{Q}_\nu^n} \bigl\langle
\Gamma,W/ \hat{S}^R \hat{T}^Q \bigr\rangle \Bigr) -
h_{\mathcal{F}_{\mu\nu}^\omega}(\Gamma) \\
&&\qquad\leq2 \bigl\{ m^{-1} |\mathcal{I}| +
n^{-1} |\mathcal{J}| \bigr\}
+ \mathbb{E} \Bigl( \max_{(Q,R) \in\mathcal{Q}_\mu^m
\times\mathcal{Q}_\nu^n}
\delta_{UV} \Bigr) \\
&&\quad\qquad{}+ 2 K \mathbb{E} \Bigl( \max
_{(Q,R) \in\mathcal{Q}_\mu^m
\times\mathcal{Q}_\nu^n} \delta_U + \delta_V \Bigr).
\end{eqnarray*}
\end{lemma}

\begin{pf}
Let $\mathring{\sigma}$ and $\mathring{\tau}$ denote the mappings in
$\mathcal{Q}_\mu$ and $\mathcal{Q}_\nu$ that are respectively closest
in the metric $d_\mathrm{Ham}$ to $\hat{\sigma}$ and $\hat{\tau}$.
Observe that we may then expand and upper-bound the left-hand side of
the lemma statement by
\begin{eqnarray*}
&&\underbrace{ \mathbb{E} \Bigl( \max_{Q,R}
\bigl\langle \Gamma,W/\hat{S}^R\hat{T}^Q \bigr\rangle -
G_{\hat{\sigma}\hat{\tau
}}(\xi, \zeta) \Bigr) }_{\mathrm{(i)}} + \underbrace{
\mathbb{E} \Bigl( \max_{Q,R} G_{\hat
{\sigma}\hat{\tau}}(
\xi,\zeta) - \langle\Gamma,\omega/ \hat {\sigma}\hat{\tau} \rangle \Bigr)
}_\mathrm{(ii)}
\\
&&\qquad{}+ \underbrace{\mathbb{E} \Bigl( \max_{Q,R}
\langle \Gamma,\omega/\hat{\sigma}\hat{\tau} \rangle - \langle\Gamma,\omega/
\mathring{\sigma} \mathring{\tau} \rangle \Bigr) }_\mathrm{(iii)} +
\underbrace{ \mathbb{E} \Bigl( \max_{Q,R}
\langle\Gamma,\omega/\mathring{ \sigma}\mathring{\tau} \rangle \Bigr) -
h_{\mathcal{F}_{\mu\nu
}^\omega}(\Gamma)}_\mathrm{(iv)},
\end{eqnarray*}
after which we may upper-bound terms (i)--(iv) in turn as follows.

First, since $|\omega(x,y)\Gamma_{\hat{\sigma}(x)\hat{\tau}(y)}|
\leq
1$ for all $(x,y)$, it follows from their respective definitions that
$\langle\Gamma,W/\hat{S}^R\hat{T}^Q \rangle - G_{\hat{\sigma
}\hat{\tau}}(\xi,\zeta
)$ is deterministically bounded above by $|\mathcal{I}| / m +
|\mathcal
{J}| / n$. Hence, term (i) is bounded by the same quantity.

Second, observe that by definition, $G_{\hat{\sigma}\hat{\tau}}(\xi
,\zeta) -\mathbb{E} ( G_{\hat{\sigma}\hat{\tau
}}(\xi
,\zeta) \vert U,V ) \leq\delta_{UV}$. Since for fixed $\sigma
, \tau$ we have $\mathbb{E} ( G_{\sigma\tau}(\xi
,\zeta) \vert U,V ) = [ |\overline{\mathcal{I}}||\overline
{\mathcal{J}}| / (mn) ] \langle\Gamma,\omega/\sigma\tau \rangle
$, with
$|\langle\Gamma,\omega/\sigma\tau \rangle| \leq1$, it holds
deterministically
that $\mathbb{E} (G_{\hat{\sigma}\hat{\tau}}(\xi
,\zeta) \vert U,V ) -\break  \langle\Gamma,\omega/\hat{\sigma}\hat
{\tau } \rangle \leq|\mathcal{I}| / m + |\mathcal{J}| / n$. Thus
term (ii) is
bounded above by the quantity\break  $\mathbb{E} ( \max_{(Q,R) \in\mathcal{Q}_\mu^m \times\mathcal{Q}_\nu^n} \delta_{UV}
) + |\mathcal{I}| / m + |\mathcal{J}| / n $.

Third, by the second Lipschitz condition of Lemma~\ref{le:lipschitz},
we have that $\langle\Gamma,\omega/\hat{\sigma}\hat{\tau}
\rangle - \langle\Gamma,\omega/\mathring{\sigma}\mathring{\tau
} \rangle \leq2 [d_\mathrm
{Ham}(\hat{\sigma},\mathring{\sigma}) + d_\mathrm{Ham}(\hat{\tau
},\mathring{\tau}) ]$. Observe that
\[
d_\mathrm{Ham}(\hat{\sigma},\mathring{\sigma}) \leq\sum
_{a=1}^K \bigl| \bigl|\hat{\sigma}^{-1}(a)\bigr| -
\mu_a \bigr| \leq\sum_{a=1}^K \Biggl| \bigl|
\hat{\sigma}^{-1}(a)\bigr| - \frac{1}{m}\sum
_{i=1}^m 1 \bigl\{\hat{\sigma}(\xi_i )
= a \bigr\} \Biggr| \leq K \delta_U,
\]
where the second inequality holds as $\hat{S}^R \in\mathcal{Q}_\mu^m$.
By the same argument for $d_\mathrm{Ham}(\hat{\tau},\mathring{\tau})$,
we see term (iii) is bounded by $2 K \mathbb{E} (
\max_{(Q,R) \in\mathcal{Q}_\mu^m \times\mathcal{Q}_\nu^n} \delta
_U +
\delta_V )$.

To conclude, note term (iv) is deterministically upper-bounded by $0$.
\end{pf}

We may now establish the claimed upper bound on $h_{{\mathcal{F}_{\mu
\nu}^A}}(\Gamma) -
h_{\mathcal{F}_{\mu\nu}^\omega}(\Gamma)$.

\begin{pf*}{Proof of Lemma~\ref{le:upper}}
Combining the results of Lemmas \ref{le:basic}--\ref{le:final} yields
directly that, with probability at least $1 - 2e^{-2mn\varepsilon^2/(m+n)}
- 2K^{m+n} e^{-2mn\varepsilon^2}$,
\begin{eqnarray*}
h_{\mathcal{F}_{\mu\nu}^A}(\Gamma) - h_{\mathcal{F}_{\mu\nu
}^\omega}(\Gamma) &\leq& 2\varepsilon+ K\sqrt{2
\pi} \bigl\{ |\mathcal{I}|^{-1/2} + |\mathcal{J}|^{-1/2} \bigr\} +
2 \bigl\{ m^{-1} |\mathcal{I}| + n^{-1} |\mathcal{J}| \bigr\}
\\
&&{}+ f \biggl(|\mathcal{I}|+|\mathcal{J}|,\ell, 2 \pmatrix{K\cr 2} \biggr)
\\
&&{} + 2 K \biggl\{ f \biggl(|
\mathcal{I}|+1,n, \pmatrix{K\cr 2} \biggr)+ f \biggl(|\mathcal{J}|
+1,m,\pmatrix{K\cr 2} \biggr) \biggr\},
\end{eqnarray*}
where $f(p,q,r) = 4 \{ [ p\log K + r\log(q + 1) + \log2
] / (2q) \}^{1/2} $, and $\ell= \min(m-|\mathcal{I}|, n
-|\mathcal{J}|)$ as in Lemma~\ref{le:rademacher1.easy}. Letting
$\varepsilon= n^{-1/4}$, $|\mathcal{I}| = |\mathcal{J}| = n^{1/2}$,
and fixing $m = \rho n$ as assumed in the hypothesis of Lemma~\ref
{le:upper}, it follows that for $n \geq2$,
\begin{eqnarray*}
h_{\mathcal{F}_{\mu\nu}^A}(\Gamma) - h_{\mathcal{F}_{\mu\nu
}^\omega}(\Gamma)& \leq&\frac{2+ 2K(2\pi)^{1/2}+(4\sqrt{2}+8K)(2
\log K)^{1/2}}{n^{1/4}}
\\
&&{}+ \frac{4 + 12(K^2 \log(\rho n+1) + 2)^{1/2}}{n^{1/2}}
\end{eqnarray*}
with probability at least $1 - 2e^{-\sqrt{n} [2\rho/(\rho+1)]} -
2K^{(\rho+1)n} e^{-2\rho n^{3/2}}$. Thus we have established the
claimed upper bound on $h_{\mathcal{F}_{\mu\nu}^A}(\Gamma)$ in
terms of $h_{\mathcal{F}_{\mu\nu}^\omega}(\Gamma)$.
\end{pf*}

\section{Simulation study}
\label{subsec:simulations}

We now present a brief simulation study which investigates empirical
rates of convergence as model misspecification increases. We control
the degree of misspecification through a sigmoidal functional form
$f_\beta(x)\dvtx [0,1] \rightarrow[-1/2,1/2]$, parameterized by $\beta
\geq1$,
\begin{eqnarray*}
f_\beta(x) &=& Z_\beta^{-1} \biggl( \frac{x^\beta}{x^\beta+
(1-x)^\beta
}
- \frac{1}{2} \biggr),\qquad 0 \leq x \leq1;\\
  Z_\beta&=& 4 \int
_0^{1/2}\biggl | \frac{x^\beta}{x^\beta+
(1-x)^\beta} - \frac{1}{2}
\biggr| \,dx.
\end{eqnarray*}
Each $f_\beta(x)$ describes a strictly monotone increasing sigmoidal
curve on $[0,1]$, proportional to $x - 1/2$ for $\beta= 1$ and to $1\{
x > 1/2\} - 1/2$ in the limit as $\beta\rightarrow\infty$.
Normalization by $Z_\beta$ maintains constant area under $|f_\beta|$.

To explore sparse graph regimes, we introduce an additional
$n$-dependent parameter $\rho_n \in(0,1)$, and take the outer product
$f_\beta(x) f_\beta(y)$ to obtain a separable generative function
$\rho
_n \omega_\beta(x,y) = \rho_n ( f_\beta(x) f_\beta(y) + 1/2
)$. As $\beta\rightarrow\infty$, this tends to a stochastic
co-blockmodel, with two classes of equal size.

\begin{figure}

\includegraphics{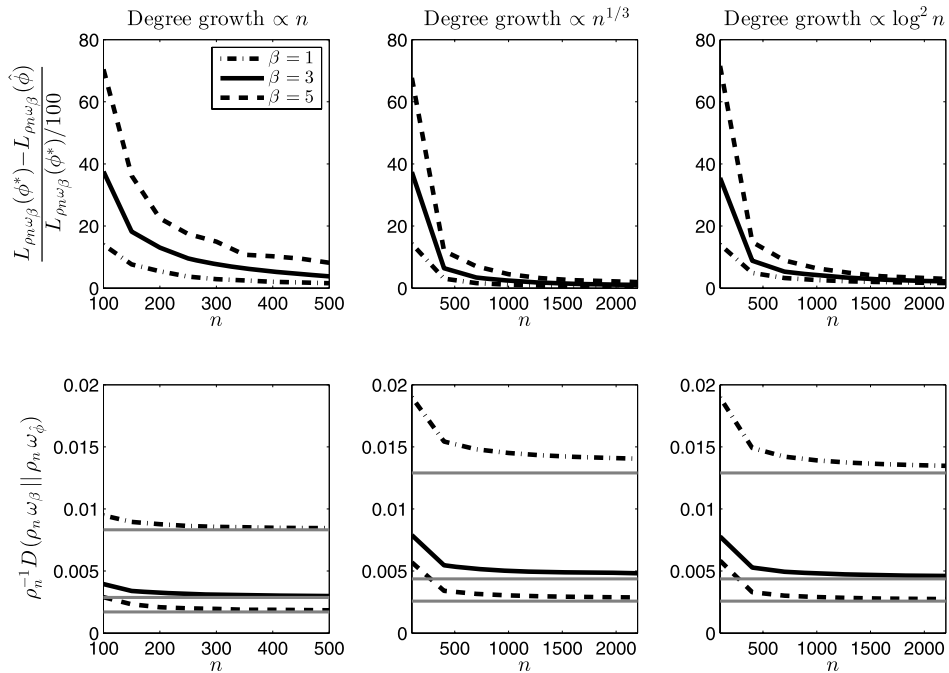}

\caption{Median performance of approximate profile
likelihood maximization according to \protect\eqref{eq:PLfunctional},
for $\rho_n \in\{ 1/2, n^{-2/3}, n^{-1}\log^2 n\}$
(left column, middle, right). Top row: percent relative excess risk,
decaying toward zero. Bottom row: Kullback--Leibler
divergence normalized by $\rho_n$,
decaying toward its asymptotic optimum in $n$ (grey horizontal
lines).}\label{fig:sims}
\end{figure}

Figure~\ref{fig:sims} shows a number of simulation results based on
this model. Specifically, for $\beta\in\{1,3,5\}$ and $\rho_n \in\{
0.5, n^{-2/3}, n^{-1}\log^2 n\}$, one thousand separable $n \times n$
binary arrays were generated from the corresponding $\rho_n \omega
_\beta(x,y)$, for network sizes ranging from 100--500 for dense graphs
(left column), and 100--2200 for sparse graphs (right columns). We see
immediately that the simulation results of Figure~\ref{fig:sims} are
qualitatively similar for all three regimes, suggesting that at least
in some cases, co-blockmodel estimators will converge despite model
misspecification in sparse as well as dense graph regimes.\looseness=1

Each of the $n \times n$ arrays described above was fitted by a
two-class co-blockmodel, whose parameters $\hat{\phi} = (\hat{\mu},
\hat
{\nu}, \hat{\theta})$ were obtained by heuristically optimizing the
profile likelihood criterion of \eqref{eq:PLfunctional} using an
algorithmic approach based on simulated annealing [\citet
{choi2012stochastic}]. Parameter values were initialized to coincide
with the optimal blockmodel approximation based on $\rho_n \omega
_\beta/\sigma^*\tau^*$, where $\sigma^*, \tau^*\dvtx [0,1] \rightarrow
\{
1,2\}$ each map the interval $[0,1/2)$ to class~$1$ and the interval
$[1/2, 1]$ to class $2$.

Lemma~\ref{le:sim1} establishes that $\phi^* = \operatorname
{argmax}_{\phi\in\Phi} L_{\rho_n\omega_\beta}(\phi)$ exists in this
setting, and that $L_{\rho_n\omega_\beta}(\phi)$ may be
straightforwardly computed for any triple $\phi= (\mu,\nu,\theta)$ of
two-class co-blockmodel parameters. Corollary~\ref{cor:sim1} then
yields a finite set containing $\phi^* = \operatorname{argmax}_{\phi
\in\Phi} L_{\rho_n\omega_\beta}(\phi)$, from which we found that
$\phi
^*$ corresponded to the blockmodel induced by $\sigma^*$ and $\tau^*$.
Thus we were able to evaluate the relative excess risk $[ L_{\rho
_n\omega_\beta}(\phi^*) - L_{\rho_n\omega_\beta}(\hat{\phi}) ]
/ L_{\rho
_n\omega_\beta}(\phi^*)$, shown as a percentage in the top row of
Figure~\ref{fig:sims}, and seen to decay toward~$0$.

The bottom row of Figure~\ref{fig:sims} shows the normalized
Kullback--Leibler divergences $\rho_n^{-1} D(\rho_n \omega_\beta
\Vert  \rho_n \omega_{\smash{\hat{\phi}}})$ decaying toward the grey
horizontal lines representing the limiting values of $D(\rho_n \omega
_\beta\Vert  \rho_n \omega_{\smash{\phi^*}})$ as $\rho_n
\rightarrow0$. These are order-one quantities, obtained through a
Taylor expansion of $D(\rho_n \omega_\beta\Vert  \rho_n \omega
_{\smash{\phi^*}})$. Smaller divergences are achieved when $\beta$ is
large, reflecting the fact that as $\beta$ increases, $\rho_n \omega
_\beta(x,y)$ becomes closer to a co-blockmodel.

Overall, we see that the simulation results shown in Figure~\ref{fig:sims} are consistent with the behavior predicted by Theorem~\ref
{th:learning} for profile likelihood maximization; qualitatively
similar results were also obtained for the least squares setting of
Theorem~\ref{th:learning} and hence are omitted for brevity.

\section{Discussion}
\label{sec:discussion}

In this article we have addressed the case of network co-clustering, in
which the inference task is to group two sets of network nodes into
classes based on their observed relations. Our results significantly
generalize known consistency results for the blockmodel and its
co-blockmodel variant: they do not require the data to be generated
(even approximately) by a co-blockmodel, and they achieve improved
rates of convergence relative to results from the graph limits
literature, through the use a Rademacher complexity bound for
$U$-statistics adapted from \citet{clemenccon2008ranking}. The assumption
of a nonparametric generative model is both more general and more
realistic, and to our knowledge Theorems \ref{th:learning} and \ref
{th:main} are the first for this regime to establish polynomial rates
of convergence.

In the work of \citet{clemenccon2008ranking}, these Rademacher
complexity results are used to derive convergence rates for learning
pairwise rankings. This setting is related to ours, but differs in some
important ways. Those authors seek a rule $r\dvtx \mathcal{X}\times
\mathcal
{X} \rightarrow\{-1,+1\}$ such that, given $X,X' \in\mathcal{X}$, $r$
indicates which has the higher rank. In this setting, $X$ and $X'$ can
be thought of as covariates describing the two objects for which a
relative ranking is desired, and $\mathcal{X}$ represents the space of
allowable covariate values. In our network setting, the nonparametric
model $\omega\dvtx [0,1]^2\rightarrow[0,1]$ is analogous to a ranking rule,
with $\mathcal{X}$ taken to be $[0,1]$. However, $X$ and $X'$ are never
observed in the data, and effectively must be imputed up to
measure-preserving transformation.

The recent work of \citet{flynn2012consistent} analyzes the consistency
of co-clustering with model misspecification, but in a rather different
setting, with the data matrix $A$ assumed to be real valued, along with
a real-valued generalization of the co-blockmodel. This generalization
utilizes discrete latent class variables $S$ and~$T$; conditioned on
$S(i)$ and $T(j)$, the distribution of $A_{ij}$ is assumed to have mean
$\theta_{S(i) T(j)}$, but may otherwise be arbitrary up to technical
conditions, and may be misspecified in the estimator. Under these
assumptions, it is shown that the latent classes can be estimated
consistently if their number is known. In the case where $A$ is binary,
the conditions of \citet{flynn2012consistent} are equivalent to
assuming a generative co-blockmodel with a known number of classes.

Finally, the very recent work of \citet{chatterjee2012matrix} derives a
simple and elegant spectral method to consistently estimate the matrix
$W$ defined in the proof of Lemma~\ref{le:upper} in Section~\ref{sec:K2Case}, that is, the mapping $\omega(x,y)$, evaluated at the
values of the latent variables $\xi_1,\ldots,\xi_m$, and $\zeta
_1,\ldots
,\zeta_n$. This implies consistency of estimation of $\omega$ in the
$L^2$ sense, and while rates of convergence are not given for general~$\omega$, they can be established for particular instances, such as
under the assumption of a generative blockmodel whose number of classes
$K$ is growing with~$n$. Our setting is distinct, in that we desire
only the best blockmodel approximation to $\omega$, and so are able to
establish $L^2$ rates of convergence that are independent
of~$\omega$.

\begin{appendix}
\section{\texorpdfstring{Proof of Theorem \lowercase{\protect\ref{th:learning}} and Lemma \lowercase{\protect\ref{lem:likelihood_support_fcns}}}
{Proof of Theorem 3.1 and Lemma 4.1}}
\label{sec:proof.th2}

To prove Theorem~\ref{th:learning}, we first denote the objective
functions of \eqref{eq:MSfunctional} and \eqref{eq:PLfunctional} by
$R_A(\phi)$ and $L_A(\phi)$, respectively. Lemma~\ref
{lem:likelihood_support_fcns}, proved below, relates $R_A(\phi) -
R_\omega
(\phi)$ and $L_A(\phi) - L_\omega(\phi)$ to the support functions
$h_{{\mathcal{F}_{\mu\nu}^A}}(\cdot)$ and $h_{{\mathcal{F}_{\mu
\nu}^\omega}}(\cdot)$, after which the result
follows directly from Theorem~\ref{th:main}.

To see this, let $\hat{\phi} \equiv(\hat{\mu},\hat{\nu},\hat
{\theta})
= \operatorname{argmin}_{\phi\in\Phi} R_A(\phi)$. For any $\phi
\in
\Phi$, we have
\begin{eqnarray*}
&&R_\omega(\hat{\phi}) - R_\omega(\phi) \\
&&\qquad =R_\omega(\hat{
\phi}) - R_A(\hat{\phi}) + R_A(\hat{\phi}) -
R_A(\phi) + R_A(\phi) - R_\omega(\phi)
\\
&&\qquad \leq R_\omega(\hat{\phi}) - R_A(\hat{\phi}) +
R_A(\phi) - R_\omega(\phi)
\\
&&\qquad \leq 2 \bigl| h_{{\mathcal{F}_{\mu\nu}^A}}( \hat{\theta}) -
 h_{{\mathcal{F}_{\mu\nu}^\omega}} ( \hat{\theta}) \bigr|
+ 2 \bigl| h_{{\mathcal{F}_{\mu\nu}^\omega}}(\theta) - h_{{\mathcal
{F}_{\mu\nu}^A}}(\theta)\bigr |,
\end{eqnarray*}
where the first inequality holds because $R_A(\hat{\phi}) - R_A(\phi)
\leq0$, and the second holds by the triangle inequality and Lemma~\ref
{lem:likelihood_support_fcns}. Applying Theorem~\ref{th:main} and
choosing $\phi$ to satisfy $R_\omega(\phi) \leq\inf_{\phi' \in
\Phi}
R_{\omega}(\phi') + n^{-1/4}$ then yields the result.

Now, assume $\phi^* = \operatorname{argmax}_{\phi\in\Phi} L_\omega
(\phi)$ exists, and set $\hat{\phi} = \operatorname{argmax}_{\phi
\in
\Phi} L_A(\phi)$. Whenever $0 < \smash{\hat{\theta}}_{ab}, {\theta
^*}_{ ab} <1$ for all $a,b=1,\ldots,K$, the second result $[ L_\omega
(\phi^*) - L_\omega(\hat{\phi}) ] / [ B(\theta^*) + B(\hat{\theta
}) ] =
\mathcal{O}_P (n^{-1/4} )$ of Theorem~\ref{th:learning} follows
similarly from
\begin{eqnarray*}
0 &\leq& L_\omega \bigl(\phi^* \bigr) - L_\omega(\hat{\phi}) \\
& =&
L_\omega \bigl(\phi^* \bigr) - L_A \bigl(\phi^* \bigr) +
L_A \bigl(\phi^* \bigr) - L_A(\hat{\phi}) +
L_A(\hat{\phi}) - L_\omega(\hat{\phi})
\\
& \leq& L_\omega \bigl(\phi^* \bigr) - L_A \bigl(\phi^*
\bigr) + L_A(\hat{\phi}) - L_\omega(\hat{\phi})
\\
& \leq& B \bigl(\theta^* \bigr)\bigl | h_{\mathcal{F}_{ \mu^* \nu^*
}^\omega
} (\Gamma_{\theta^* }) -
h_{\mathcal{F}_{ \mu^* \nu^* }^A} (\Gamma_{\theta^* })\bigr | + B(\hat {\theta})\bigl |
h_{\mathcal{F}_{\hat{\mu} \hat{\nu}}^A}(\Gamma_{\hat{\theta}}) - h_{\mathcal{F}_{\hat{\mu}\hat{\nu}}^\omega
}(
\Gamma_{\hat{\theta}}) \bigr|.
\end{eqnarray*}

\begin{pf*}{Proof of Lemma~\ref{lem:likelihood_support_fcns}}
We show the results of the lemma directly,
\begin{eqnarray*}
R_A(\phi) & =& \min_{(S,T) \in\mathcal{Q}_\mu^m \times\mathcal
{Q}_\nu^n} \frac
{1}{mn} \sum
_{i=1}^m \sum
_{j=1}^n \llvert\theta_{S(i) T(j)} -
A_{ij} \rrvert^2
\\
& = &\min_{F \in{\mathcal{F}_{\mu\nu}^A}} \Biggl\{ \sum_{a=1}^K
\sum_{b=1}^K -2 F_{ab}
\theta_{ab} + \mu_a \nu_b
\theta_{ab}^2 \Biggr\} + \frac{1}{mn} \sum
_{i=1}^m \sum_{j=1}^n
A_{ij}^2
\\
& = &\Biggl\{ - 2 h_{{\mathcal{F}_{\mu\nu}^A}}(\theta) + \sum_{a=1}^K
\sum_{b=1}^K \mu_a \nu
_b \theta_{ab}^2 \Biggr\} + \frac{1}{mn}
\sum_{i=1}^m \sum
_{j=1}^n A_{ij}^2,
\end{eqnarray*}
where the second line follows from the definition of ${\mathcal
{F}_{\mu\nu}^A}$, and the
last line from that of $h_{{\mathcal{F}_{\mu\nu}^A}}$. Letting
$(\sigma,\tau)$ satisfy
$\sigma(x) = F^{\smash{-1}}_{\smash{\mu(\pi_1(x))}}$ and $\tau(y) =
F^{\smash{-1}}_{\smash{\nu(\pi_2(y))}}$,
\begin{eqnarray*}
R_\omega(\phi) & =& \inf_{\pi_1,\pi_2 \in\mathcal{P}} \int
_{[0,1]^2} \bigl\llvert\omega \bigl(\pi_1(x),
\pi_2(y) \bigr) - \omega_\phi(x,y) \bigr\rrvert
^2 \,dx \,dy
\\
& = &\inf_{(\sigma,\tau) \in\mathcal{Q}_\mu\times\mathcal{Q}_\nu
} \sum_{a=1}^K
\sum_{b=1}^K \int_{\sigma^{-1}(a)\times\tau^{-1}(b)}
\bigl\llvert\omega(x,y)- \theta_{ab} \bigr\rrvert^2 \,dx \,dy
\\
& =& \inf_{F \in{\mathcal{F}_{\mu\nu}^\omega}} \Biggl\{ \sum_{a=1}^K
\sum_{b=1}^K -2F_{ab} \theta
_{ab} + \mu_a\nu_b \theta_{ab}^2
\Biggr\} + \int_{[0,1]^2} \omega(x,y)^2 \,dx \,dy
\\
& =& \Biggl\{ - 2 h_{{\mathcal{F}_{\mu\nu}^\omega}}(\theta) + \sum_{a=1}^K
\sum_{b=1}^K \mu_a \nu
_b \theta_{ab}^2 \Biggr\} + \int
_{[0,1]^2} \omega(x,y)^2 \,dx \,dy.
\end{eqnarray*}

Following similar steps, we show the second result as follows:
\begin{eqnarray*}
L_A(\phi) & = &\max_{(S,T) \in\mathcal{Q}_\mu^m \times\mathcal
{Q}_\nu
^n} \frac{1}{mn} \sum
_{i=1}^m \sum
_{j=1}^n \bigl\{ A_{ij} \log( \theta
_{S(i)T(j)}) \\
&&\hspace*{110pt}{}+ (1-A_{ij}) \log(1-\theta_{S(i)T(j)}) \bigr\}
\\
& = &\max_{F \in{\mathcal{F}_{\mu\nu}^A}} \sum_{a=1}^K
\sum_{b=1}^K \biggl\{ F_{ab}
\log \biggl( \frac{\theta_{ab}}{1-\theta_{ab}} \biggr) + \mu_a \nu_b
\log(1-\theta_{ab}) \biggr\}
\\
& =& B(\theta) h_{{\mathcal{F}_{\mu\nu}^A}}(\Gamma_\theta) + \sum
_{a=1}^K \sum_{b=1}^K
\mu_a \nu_b \log(1-\theta_{ab}),
\end{eqnarray*}
since $\max_{F \in{\mathcal{F}_{\mu\nu}^A}} \sum_{a,b} F_{ab}
B(\theta) (\Gamma
_\theta)_{ab} = B(\theta) h_{{\mathcal{F}_{\mu\nu}^A}}(\Gamma
_\theta)$, and similarly
\begin{eqnarray*}
L_\omega(\phi) & =& \sup_{\pi_1,\pi_2 \in\mathcal{P}} \int
_{[0,1]^2} \bigl\{ \omega \bigl(\pi_1(x),
\pi_2(y) \bigr) \log\omega_\phi(x,y)
\\
&&\hspace*{62pt}{}  + \bigl[1-\omega \bigl(\pi_1(x), \pi_2(y)
\bigr) \bigr] \log \bigl(1-\omega_\phi(x,y) \bigr) \bigr\} \,dx \,dy
\\
& = &\sup_{(\sigma,\tau) \in\mathcal{Q}_\mu\times\mathcal{Q}_\nu
} \sum_{a=1}^K
\sum_{b=1}^K \int_{\sigma^{-1}(a)\times\tau^{-1}(b)}
\bigl\{ \omega(x,y) \log\theta_{ab}\\
&&\hspace*{154pt}{}+ \bigl(1-\omega(x,y) \bigr) \log(1-\theta_{ab}) \bigr
\} \,dx \,dy
\\
& =& \sup_{F \in{\mathcal{F}_{\mu\nu}^\omega}} \sum_{a=1}^K
\sum_{b=1}^K \biggl\{ F_{ab}
\log \biggl( \frac{\theta_{ab}}{1-\theta_{ab}} \biggr) + \mu_a \nu_b
\log(1-\theta_{ab}) \biggr\}
\\
& = &B(\theta) h_{{\mathcal{F}_{\mu\nu}^\omega}}( \Gamma_\theta) + \sum
_{a=1}^K \sum_{b=1}^K
\mu_a \nu_b \log(1-\theta_{ab}).
\end{eqnarray*}
\upqed\end{pf*}

\section{\texorpdfstring{Auxiliary proofs for Theorem \lowercase{\protect\ref{th:main}}}
{Auxiliary proofs for Theorem 4.1}}
\label{sec:proof.th3}

Below we provide proofs of all supporting lemmas for Theorem~\ref
{th:main}, and state and prove the covering argument used to establish
the theorem:
\begin{longlist}[(1)]

\item[(1)] First, in Sections~\ref{le:lowerpf}--\ref{le:alonpf} below,
we prove auxiliary Lemmas \ref{le:lower}, \ref{le:basic} and~\ref
{le:alon} as stated in Section~\ref{subsec:sketch}.

\item[(2)] Then, in Section~\ref{sec:QU and QV}, we generalize the
definitions of $\mathcal{Q}_U$ and $\mathcal{Q}_V$, given in
Section~\ref{sec:K2Case} for $K = 2$, to arbitrary $K$; this induces
generalizations of the quantities $\delta_{U}, \delta_{V}$ and
$\delta
_{UV}$ in the natural way.

\item[(3)] Then, in Sections~\ref{le:rademacher1.easypf} and \ref
{le:rademacher2.easypf}, we prove Lemmas \ref{le:rademacher1.easy}
and \ref{le:rademacher2.easy}, which depend on $(\mathcal{Q}_U,
\mathcal
{Q}_V, \delta_{U}, \delta_{V}, \delta_{UV})$ as defined for
arbitrary $K$.

\item[(4)] Finally, in Section~\ref{sec:epsCovering}, we extend the
pointwise convergence result of Proposition~\ref{prop:hyperplane} by
way of a covering argument for all $\Gamma\in[-1,1]^{K \times K}$.
\end{longlist}

\subsection{\texorpdfstring{Proof of Lemma \protect\ref{le:lower}}
{Proof of Lemma 5.2}}
\label{le:lowerpf}

For fixed $\Gamma$, let $(\sigma^*,\tau^*) \in\mathcal{Q}_\mu
\times
\mathcal{Q}_\nu$ satisfy
%
\begin{equation}
\label{eq:sketch0} \bigl\langle\Gamma,\omega/\sigma^*\tau^* \bigr\rangle >
h_{{\mathcal{F}_{\mu\nu}^\omega}}(\Gamma) - \frac{1}{n^{1/4}},
\end{equation}
so that $\omega/\sigma^* \tau^*$ is within $n^{-1/4}$ of the supporting
hyperplane. Define
\[
S^*(i) = \sigma^*(\xi_i), \qquad T^*(j) = \tau^*(\zeta_j);\qquad
i=1,\ldots,m, j=1,\ldots,n.
\]
By the arguments of Lemma~\ref{le:basic} as proved in Section~\ref{le:basicpf} below, applying McDiarmid's inequality with the Lipschitz
conditions of Lemma~\ref{le:lipschitz} yields
%
\begin{equation}
\label{eq:sketch1}\qquad \mathbb{P} \bigl(\bigl|\bigl\langle\Gamma,A/S^*T^* \bigr\rangle -
\bigl\langle\Gamma,\omega/ \sigma^*\tau^* \bigr\rangle\bigr| \geq2 \varepsilon \bigr)
\leq2e^{-2mn\varepsilon^2/(m+n)} + 2e^{-2mn\varepsilon^2}.
\end{equation}
While $(S^*,T^*)$ many not be in $\mathcal{Q}_\mu^m \times\mathcal
{Q}_\nu^n$, a Chernoff bound implies that
\[
\mathbb{P} \biggl( \biggl\llvert\frac{S^{*-1}(a)}{m} - \mu_a \biggr
\rrvert\geq\varepsilon \biggr) \leq2e^{-2m\varepsilon^2},\qquad a=1,\ldots,K.
\]
The analogous bound also holds for $|T^{*-1}(b)/n - \nu_b|$. Applying
these results in conjunction with a union bound yields
\[
\mathbb{P} \biggl(\max_{1 \leq a,b \leq K} \biggl\{ \biggl\llvert
\frac{S^{*-1}(a)}{m}-\mu_a \biggr\rrvert+ \biggl\llvert
\frac
{T^{*-1}(b)}{n} - \nu_b \biggr\rrvert \biggr\} \geq2\varepsilon
\biggr) \leq K \bigl(2e^{-2m \varepsilon^2} + 2e^{-2n \varepsilon^2} \bigr).
\]
Therefore, with probability at least $1-K(2e^{-2m \varepsilon^2} + 2e^{-2n
\varepsilon^2})$, there exists a pair $(\mathring{S},\mathring{T}) \in
\mathcal{Q}_\mu^m \times\mathcal{Q}_\nu^n$ such that
\[
\frac{1}{m}d_{\mathrm{Ham}} \bigl(S^*,\mathring{S} \bigr) +
\frac{1}{n}d_{\mathrm
{Ham}} \bigl(T^*,\mathring{T} \bigr) \leq2K\varepsilon,
\]
which by the first condition of Lemma~\ref{le:lipschitz} implies that
%
\begin{equation}
\label{eq:sketch2} \bigl|\langle\Gamma,A/\mathring{S}\mathring{T} \rangle - \bigl\langle
\Gamma,A/S^*T^* \bigr\rangle\bigr| \leq4K\varepsilon.
\end{equation}

Recalling that $h_{{\mathcal{F}_{\mu\nu}^A}} = \max_{(S,T) \in
\mathcal{Q}_\mu^m \times
\mathcal{Q}_\nu^n} \langle\Gamma,A/ST \rangle$, we have that
\[
h_{{\mathcal{F}_{\mu\nu}^A}}(\Gamma)  \geq\langle\Gamma,A/\mathring{S}\mathring{T}
\rangle,
\]
following which \eqref{eq:sketch2}, \eqref{eq:sketch1}
and \eqref{eq:sketch0} in turn imply that with probability
at least $1 - 2e^{-2mn\varepsilon^2/(m+n)} - 2e^{-2mn\varepsilon^2}
- K(2e^{-2m \varepsilon^2} + 2e^{-2n \varepsilon^2})$, we have
\begin{eqnarray*}
h_{{\mathcal{F}_{\mu\nu}^A}}(\Gamma) & \geq&\bigl\langle\Gamma,A/S^*T^* \bigr\rangle - 4K
\varepsilon
\\
& \geq&\bigl\langle\Gamma,\omega/\sigma^*\tau^* \bigr\rangle - (4K+2)\varepsilon
\\
& \geq& h_{{\mathcal{F}_{\mu\nu}^\omega}}(\Gamma) - n^{-1/4} - (4K+2)\varepsilon.
\end{eqnarray*}
Now letting $m = \rho n$ as in the statement of the lemma,
and setting $\varepsilon= n^{-1/4}$, we see that with probability
at least $1 - 2e^{-\sqrt{n} [2\rho/(\rho+1)]} [1+o(1) ]$,
\[
h_{{\mathcal{F}_{\mu\nu}^A}}(\Gamma)  \geq h_{{\mathcal{F}_{\mu
\nu}^\omega}}(\Gamma) - \frac{4K+3}{n^{1/4}},
\]
providing the necessary lower bound on $h_{{\mathcal{F}_{\mu\nu
}^A}}(\Gamma)$ in terms of
$h_{{\mathcal{F}_{\mu\nu}^\omega}}(\Gamma)$.

\subsection{\texorpdfstring{Proof of Lemma \protect\ref{le:basic}}
{Proof of Lemma 5.4}}
\label{le:basicpf}

Recalling the definitions of $h_{{\mathcal{F}_{\mu\nu}^A}}$ and
$h_{{\mathcal{F}_{\mu\nu}^W}}$,
%
\begin{eqnarray}
\nonumber
&&\mathbb{P} \bigl(\bigl|h_{{\mathcal{F}_{\mu\nu}^A}}(\Gamma) -
h_{{\mathcal{F}_{\mu\nu}^W}}(
\Gamma)\bigr| \geq\varepsilon \bigr)
\\
\nonumber
&&\qquad = \mathbb{P} \Bigl( \Bigl\llvert\max_{(S,T) \in\mathcal{Q}_\mu
^m \times\mathcal{Q}_\nu^n}
\langle\Gamma,A/ST \rangle - \max_{(S,T) \in\mathcal
{Q}_\mu^m \times\mathcal{Q}_\nu^n} \langle\Gamma,W/ST
\rangle \Bigr\rrvert\geq \varepsilon \Bigr)
\\
\nonumber
&&\qquad \leq\mathbb{P} \Bigl(\max_{(S,T) \in\mathcal{Q}_\mu^m
\times\mathcal{Q}_\nu^n} \bigl\llvert
\langle\Gamma,A/ST \rangle - \langle\Gamma,W/ST \rangle \bigr\rrvert\geq\varepsilon
\Bigr)
\\
\label{eq:hAWbnd} &&\qquad \leq\sum_{(S,T) \in\mathcal{Q}_\mu^m \times
\mathcal{Q}_\nu^n} \mathbb{P} \bigl(
\bigl\llvert\langle\Gamma,A/ST \rangle - \langle\Gamma,W/ST \rangle \bigr\rrvert
\geq\varepsilon \bigr)
\\
\label{eq:hAWbndE} &&\qquad = \sum_{(S,T) \in\mathcal{Q}_\mu^m \times
\mathcal{Q}_\nu^n} \mathbb{P} \bigl( \bigl|
\langle\Gamma,A/ST \rangle - \mathbb{E} \bigl( \langle \Gamma,A/ST
\rangle \bigr)\bigr |\geq\varepsilon \bigr),
\end{eqnarray}
where \eqref{eq:hAWbnd} follows by a union bound, and \eqref
{eq:hAWbndE} by considering $\langle\Gamma,A/ST \rangle$ as a
function of the
$mn$ independent random variables $\{A_{ij}\}$, which shows that
$\mathbb{E} ( \langle\Gamma,A/ST \rangle ) =
\langle\Gamma,W/ST \rangle$ for each $(S,T)$, as $W_{ij} = \omega
(\xi_i,\zeta_j) =
\mathbb{E} (A_{ij} )$.

Next, recall the final Lipschitz condition of Lemma~\ref{le:lipschitz},
which states that $\llvert\langle\Gamma,A/ST \rangle - \langle
\Gamma,A'/ST \rangle\rrvert
\leq1/(mn)$ if $A$ and $A'$ differ by a single entry. Thus we may
apply McDiarmid's inequality to bound each term in \eqref{eq:hAWbndE},
and since $|\mathcal{Q}_\mu^m| \leq K^m$ and $|\mathcal{Q}_\nu^n|
\leq
K^n$, we obtain after summing that
\[
\mathbb{P} \bigl(\bigl|h_{{\mathcal{F}_{\mu\nu}^A}}(\Gamma)
 - h_{{\mathcal{F}_{\mu\nu}^W}}(\Gamma)\bigr| \geq
\varepsilon \bigr) \leq K^{m+n} \cdot2e^{-2mn\varepsilon^2}.
\]

Now consider $h_{{\mathcal{F}_{\mu\nu}^W}}(\Gamma) = \max_{(S,T)
\in\mathcal{Q}_\mu^m
\times\mathcal{Q}_\nu^n} \langle\Gamma,W/ST \rangle$ as a
function of the $m+n$
independent random variables $\xi_1,\ldots,\xi_m$ and $\zeta
_1,\ldots
,\zeta_n$. Changing a single component of $\xi$ or $\zeta$ affects only
a single row or column of $W$, respectively, and thus alters $\langle
\Gamma,W/ST \rangle$ and hence $h_{\mathcal{F}_{\mu\nu}^W}$ by at
most $1/m$ or
$1/n$. It therefore follows directly from McDiarmid's inequality that
\[
\mathbb{P} \bigl( \bigl|h_{\mathcal{F}_{\mu\nu}^W}(\Gamma) - \mathbb{E}
h_{\mathcal{F}_{\mu\nu}^W}(\Gamma)\bigr| \geq\varepsilon \bigr) \leq 2e^{-2mn\varepsilon^2/(m+n)}.
\]
Combining these inequalities via a union bound yields the statement of
the lemma, since by the triangle inequality we must have $|h_{\smash
{{\mathcal{F}_{\mu\nu}^A}}}(\Gamma) - h_{\smash{{\mathcal{F}_{\mu
\nu}^W}}}(\Gamma)| \geq\varepsilon$ or
$|h_{\mathcal{F}_{\mu\nu}^W}(\Gamma) - \mathbb{E}
h_{\mathcal{F}_{\mu\nu}^W}(\Gamma)| \geq
\varepsilon$ in order that $| h_{{\mathcal{F}_{\mu\nu}^A}}(\Gamma) -
\mathbb
{E} h_{\mathcal{F}_{\mu\nu}^W}(\Gamma) | \geq2\varepsilon$.

\subsection{\texorpdfstring{Proof of Lemma \protect\ref{le:alon}}
{Proof of Lemma 5.5}}
\label{le:alonpf}

Recall from the statement of the lemma that $\mathcal{I}$ and
$\mathcal
{J}$ denote sets of deterministic size whose elements are sampled
without replacement from $1,\ldots,m$ and $1,\ldots,n$, respectively.
We adopt the notation that $\mathbb{E}_{\mathcal{I}}$
denotes an expectation taken over $\mathcal{I}$, with all other random
variables held constant, and define $\mathbb
{E}_{\mathcal
{J}}$ and $\mathbb{E}_{\mathcal{I}\mathcal{J}}$ in the
same manner.

To prove the lemma, it suffices to show that for all $W,T,S$,
%
\begin{eqnarray}
\label{eq:alon1} \mathbb{E}_{\mathcal{J}} \bigl( \bigl\langle
\Gamma,W/ \hat{S}^T T \bigr\rangle \bigr) & \geq&\bigl\langle\Gamma
,W/S^TT \bigr\rangle -K\sqrt{ 2 \pi/ |\mathcal{J}| },
\\
\label{eq:alon2} \mathbb{E}_{\mathcal{I}} \bigl(\bigl\langle
\Gamma, W/S \hat{T}^S \bigr\rangle \bigr) & \geq&\bigl\langle\Gamma,
W/ST^S \bigr\rangle - K\sqrt{ 2 \pi/ |\mathcal{I}| },
\end{eqnarray}
where $\hat{S}^T$ and $\hat{T}^S$ are respectively defined in \eqref
{eq:hatS} and \eqref{eq:hatT}, and
\[
S^T = \mathop{\operatorname{argmax}}_{S \in\mathcal{Q}_\mu^m} \langle\Gamma,W/ST
\rangle, \qquad T^S = \mathop{\operatorname{argmax}}_{T \in\mathcal{Q}_\nu^n} \langle\Gamma
,W/ST \rangle.
\]
This is because \eqref{eq:alon1} and \eqref{eq:alon2} imply that for
all $(U,V) \in\mathcal{Q}_\mu^m\times\mathcal{Q}_{\nu}^n$,
\begin{eqnarray*}
\langle\Gamma, W/UV \rangle &\leq&\bigl\langle\Gamma, W/UT^U \bigr
\rangle
\\
& \leq&\mathbb{E}_{\mathcal{I}} \bigl(\bigl\langle\Gamma, W/U
\hat{T}^U \bigr\rangle \bigr) + K\sqrt{ 2\pi/ |\mathcal{I}| }
\\
& \leq&\mathbb{E}_{\mathcal{I}} \bigl(\bigl\langle\Gamma,
W/S^{\hat{T}^U}\hat{T}^U \bigr\rangle \bigr) + K\sqrt{ 2\pi/ |
\mathcal{I}| }
\\
& \leq&\mathbb{E}_{\mathcal{I}}\mathbb
{E}_{\mathcal{J}} \bigl(\bigl\langle\Gamma, W/\hat{S}^{\hat
{T}^U}\hat
{T}^U \bigr\rangle \bigr) + K\sqrt{ 2\pi/ |\mathcal{I}| } + K\sqrt{
2\pi/ | \mathcal{J}| }
\\
& \leq&\mathbb{E}_{\mathcal{I}\mathcal{J}} \Bigl(\max_{(Q,R)
 \in\mathcal{Q}_\mu^m\times\mathcal{Q}_{\nu}^n}
\bigl\langle\Gamma, W/\hat{S}^{R}\hat{T}^Q \bigr\rangle
\Bigr)\\
&&{} + K\sqrt {2\pi} \bigl( |\mathcal{I}|^{-{1}/{2}} + |
\mathcal{J}|^{-{1}/{2}} \bigr).
\end{eqnarray*}
Recalling the definition of $h_{\mathcal{F}_{\mu\nu}^W}(\Gamma)$,
and noting that the
right-hand side above is deterministic for fixed $W$, with no
dependence on $U$ or $V$, we may write
\begin{eqnarray*}
h_{\mathcal{F}_{\mu\nu}^W}(\Gamma) & =& \max_{(U,V) \in\mathcal
{Q}_\mu^m \times\mathcal
{Q}_\nu^n} \langle\Gamma,W/UV
\rangle
\\
& \leq&\mathbb{E}_{\mathcal{I}\mathcal{J}} \Bigl(\max_{(Q,R) \in\mathcal{Q}_\mu^m\times\mathcal{Q}_{\nu}^n}
\bigl\langle\Gamma, W/\hat{S}^{R}\hat{T}^Q \bigr\rangle
\Bigr)\\
&&{} + K\sqrt {2\pi} \bigl( |\mathcal{I}|^{-{1}/{2}} + |
\mathcal{J}|^{-{1}/{2}} \bigr).
\end{eqnarray*}
Taking expectations on both sides over $W$ gives the statement of the lemma.

We now establish \eqref{eq:alon1}, noting that \eqref{eq:alon2} will
follow by parallel arguments. For fixed $W$ and $T$, define for any $a
= 1, \ldots, K$ the difference
\[
\Delta_i^a = \frac{1}{|\mathcal{J}|} \sum
_{j\in\mathcal{J}} W_{ij} \Gamma_{aT(j)} -
\frac{1}{n} \sum_{j=1}^n
W_{ij} \Gamma_{aT(j)}.
\]
It follows that $\mathbb{E}_{\mathcal{J}} ( \Delta
_i^a ) = 0$, and by a Chernoff bound,
\[
\mathbb{P} \bigl( \bigl|\Delta_i^a\bigr| \geq t \bigr)
\leq2e^{-2t^2|\mathcal{J}|}.
\]
As $|\Delta_i^a|$ is nonnegative, the identity $\mathbb
{E}(X) = \int_0^\infty\mathbb{P}(X \geq t) \,dt$ for $X$ taking only
nonnegative values can be used to bound its expectation according to
\[
\mathbb{E}_\mathcal{J} \bigl(\bigl|\Delta_i^a\bigr|
\bigr) \leq\sqrt{\pi/\bigl(2|\mathcal{J}|\bigr)},
\]
which implies
%
\begin{equation}
\label{eq:approx_delta} \mathbb{E}_{\mathcal{J}} \Bigl(\max
_{1 \leq a \leq K} \bigl|\Delta_i^a\bigr| \Bigr) \leq K
\sqrt{\pi/\bigl(2|\mathcal{J}|\bigr)}.
\end{equation}

For fixed $W$ and $\mathcal{J}$, define the function
\[
f_W(S,T) = \frac{1}{m|\mathcal{J}|} \sum_{i=1}^m
\sum_{j \in\mathcal
{J}} W_{ij} \Gamma_{S(i)T(j)},
\]
and for fixed $W$ and $T$, let
%
\begin{equation}
\label{eq:deltaDefn} \Delta= \max_{S\in\mathcal{Q}_\mu^m} \bigl\llvert
f_W(S,T) - \langle\Gamma,W/ST \rangle \bigr\rrvert.
\end{equation}
From the definition of $\Delta$ it follows that
\begin{eqnarray*}
\Delta&=& \max_{S\in\mathcal{Q}_\mu^m} \Biggl\{ \frac{1}{m} \Biggl\llvert
\sum_{i=1}^m \Biggl( \frac{1}{|\mathcal{J}|}
\sum_{j\in\mathcal{J}} W_{ij} \Gamma_{S(i)T(j)} -
\frac{1}{n} \sum_{j=1}^n
W_{ij} \Gamma_{S(i)T(j)} \Biggr) \Biggr\rrvert \Biggr\}
\\
& \leq&\frac{1}{m} \Biggl\llvert\sum_{i=1}^m
\max_{1 \leq a \leq K} \Biggl\{ \frac{1}{|\mathcal{J}|} \sum
_{j\in\mathcal{J}} W_{ij} \Gamma_{aT(j)} -
\frac{1}{n} \sum_{j=1}^n
W_{ij} \Gamma_{aT(j)} \Biggr\} \Biggr\rrvert
\\
& = &\frac{1}{m} \Biggl| \sum_{i=1}^m  \max_{1 \leq a \leq K} \bigl\{ \Delta_i^a \bigr
\} \Biggr\rrvert\leq\frac{1}{m} \sum_{i=1}^m
\max_{1 \leq a \leq K} \bigl\llvert\Delta_i^a
\bigr\rrvert.
\end{eqnarray*}
Taking expectations of both sides over $\mathcal{J}$ and
substituting \eqref{eq:approx_delta} yields
%
\begin{equation}
\label{eq:approx_Delta} \mathbb{E}_{\mathcal{J}} ( \Delta) \leq
\frac{1}{m} \sum_{i=1}^m
\mathbb{E}_{\mathcal{J}} \Bigl(\max_{1 \leq a \leq
K} \bigl|
\Delta_i^a \bigr| \Bigr) \leq K \sqrt{\pi/\bigl(2|\mathcal{J}|\bigr)}.
\end{equation}

Finally, to show \eqref{eq:alon1}, observe that since $\hat{S}^T$
from \eqref{eq:hatS} maximizes $f_W(\cdot,T)$, and $S^T$ as defined
above maximizes $\langle\Gamma, W/\cdot T \rangle$, we have from
\eqref
{eq:deltaDefn} that
\begin{eqnarray*}
0 & \leq&\bigl\langle\Gamma, W/S^TT \bigr\rangle - \bigl\langle
\Gamma, W/\hat {S}^T T \bigr\rangle
\\
& \leq&\bigl\langle\Gamma, W/S^TT \bigr\rangle - f_W
\bigl(S^T,T \bigr) + f_W \bigl(\hat{S}^T,T
\bigr) - \bigl\langle\Gamma, W/ \hat {S}^T T \bigr\rangle \leq2
\Delta,
\end{eqnarray*}
and so $\langle\Gamma, W/\hat{S}^T T \rangle \geq\langle\Gamma,
W/S^TT \rangle - 2 \Delta$.
Taking expectations of both sides of this expression over $\mathcal
{J}$, and then substituting \eqref{eq:approx_Delta}, yields the inequality
\[
\mathbb{E}_{\mathcal{J}} \bigl( \bigl\langle\Gamma, W/
\hat{S}^T T \bigr\rangle \bigr) \geq\bigl\langle\Gamma,
W/S^TT \bigr\rangle - 2 K \sqrt{ \pi/\bigl(2|\mathcal{J}|\bigr)},
\]
which is the statement of \eqref{eq:alon1}. That of \eqref{eq:alon2}
follows by parallel arguments.

\subsection{\texorpdfstring{Definition of $\mathcal{Q}_U$ and $\mathcal{Q}_V$ for arbitrary $K$}
{Definition of Q U and Q V for arbitrary K}}
\label{sec:QU and QV}

In order to redefine $\mathcal{Q}_U$ and $\mathcal{Q}_V$ to accommodate
arbitrary $K$, we first redefine the mappings $U$ and $V$. Given $\zeta
_\mathcal{J} = \{\zeta_j\dvtx j \in\mathcal{J}\}$ and an assignment
function $R\dvtx \{1,\ldots,n\} \rightarrow\{1,\ldots,K\}$, define the
mapping $U\dvtx [0,1]\rightarrow\mathbb{R}^K$ by
\[
U_a(x) = \sum_{j\in\mathcal{J}} \omega(x,
\zeta_j) \Gamma_{aR(j)}, \qquad x \in[0,1], a=1,\ldots,K.
\]
Analogously, given $\xi_{\mathcal{I}}$ and $Q$, define
$V\dvtx [0,1]\rightarrow\mathbb{R}^K$ by
\[
V_a(y) = \sum_{i \in\mathcal{I}} \omega(
\xi_i,y) \Gamma_{Q(i) a},\qquad  y \in[0,1], a=1,\ldots,K.
\]

Given $a,b \in\{1,\ldots,K\}$ and the mapping $U$, define the relation
$\succeq^{U,a,b}$ by
\begin{eqnarray*}
x_1 &\succeq&^{U,a,b} x_2 \\
&\equiv&\cases{
U_a(x_1) - U_b(x_2) >
U_a(x_2) - U_b(x_1), & \quad $\mbox{or}$ \vspace*{2pt}
\cr
U_a(x_1) -
U_b(x_2) = U_a(x_2) -
U_b(x_1), &\quad  $\mbox{if $(a - b) (x_1 -
x_2) \geq0$.}$}
\end{eqnarray*}
Informally, $x_1 \succeq^{U,a,b} x_2$ implies that, given the choice of
assigning either $x_1$ or $x_2$ to group $a$, with the other relegated
to group $b$, $x_1$ is at least as attractive as $x_2$. The latter
tie-breaker condition results in a symmetric definition: if $x_1
\succeq
^{U,a,b} x_2$, then $x_2 \succeq^{U,b,a} x_1$. We define $\succ
^{U,a,b}$ analogously to $\succeq^{U,a,b}$, except that the inequality
$(a-b)(x_1-x_2) > 0$ is strict.

Let $\mathcal{S}$ denote the set of symmetric matrices in $[0,1]^{K
\times K}$. Given $t \in\mathcal{S}$ and the mapping $U$, we define
the function $\sigma_t\dvtx [0,1]\rightarrow\{1,\ldots,K\}$ as the mapping
which satisfies the following:
\[
\sigma_t^{-1}(a) = \bigl\{x\dvtx x \succeq^{U,a,b}
t_{ab}\ \forall b > a, x \succ^{U,a,b} t_{ab}\ \forall b < a \bigr\},\qquad a=1,\ldots,K,
\]
with the convention that $\sigma_t$ is undefined whenever the above
rule does not map all of $[0,1]$ to $\{1,\ldots,K\}$.

We define the function class $\mathcal{Q}_U$ as follows:
\[
\mathcal{Q}_U = \{\sigma_t \dvtx t \in\mathcal{S} \mbox{
and $\sigma_t$ is defined} \}.
\]
Given $t \in\mathcal{S}$ and the mapping $V$ as defined above, we
define $\succ^{V,a,b}, \tau_t$ and $\mathcal{Q}_V$ analogously. We then
have the following.

\begin{lemma} \label{le:hatS}
Given $U$ induced by $\zeta_\mathcal{J}$ and $R$, and given $W$ induced
by $\xi$ and~$\zeta$, define $\hat{S}^R$ by \eqref{eq:hatS}. Then there
exists $\hat{\sigma} \in\mathcal{Q}_U$ such that
\[
\hat{S}^R(i) = \hat{\sigma}(\xi_i),\qquad i=1,\ldots,m.
\]
Likewise, given $V$ induced by $\xi_\mathcal{I}$ and $Q$, and given $W$
induced by $\xi$ and $\zeta$, define $\hat{T}^Q$ by \eqref{eq:hatT}.
Then there exists $\hat{\tau} \in\mathcal{Q}_V$ such that
\[
\hat{T}^Q(j) = \hat{\tau}(\zeta_j),\qquad j=1,\ldots,n.
\]
\end{lemma}

\begin{pf}
Let $\hat{S}^R$ be chosen lexicographically from the set of all
maximizers of~\eqref{eq:hatS}, where $S$ lexicographically precedes
$S'$ if and only if $S(i_1), \ldots, S(i_m)$ lexicographically precedes
$S'(i_1),\ldots,S(i_m)$, where $i_1,\ldots,i_m$ are in order of
increasing $\xi_{i_1},\ldots,\xi_{i_m}$.

Since $\hat{S}^R$ maximizes \eqref{eq:hatS}, it holds for all $i,j
=1,\ldots,m$ that
\[
U_{\hat{S}^R(i)}(\xi_i) + U_{\hat{S}^R(j)}(\xi_j)
\geq U_{\hat
{S}^R(i)}(\xi_j) + U_{\hat{S}^R(j)}(
\xi_i);
\]
otherwise switching labels for $i$ and $j$ would increase the value of
the objective function. As $\hat{S}^R$ is chosen lexicographically, for
any $i,j$ such that
\[
U_{\hat{S}^R(i)}(\xi_i) + U_{\hat{S}^R(j)}(\xi_j) =
U_{\hat{S}^R(i)}(\xi_j) + U_{\hat{S}^R(j)}(\xi_i),
\]
it holds that $ (\hat{S}^R(i) - \hat{S}^R(j) ) (\xi_i - \xi_j)
\geq0$, with equality if and only if $\xi_i=\xi_j$. Otherwise,
switching labels would improve the lexicographic ordering.

Since $\xi_i \neq\xi_j$ for $i\neq j$ except on a set of measure zero,
it follows that
\[
\bigl(\hat{S}^R \bigr)^{-1}(a) \succ^{U,a,b} \bigl(
\hat{S}^R \bigr)^{-1}(b),\qquad  a,b = 1,\ldots,K, a\neq b,
\]
where we have let $ (\hat{S}^R )^{-1}(a)$ denote $\{\xi_i\dvtx \hat
{S}^R(\xi_i) =a \}$. As a result, for each $a$ and $b$ we may choose
$t_{ab}=t_{ba} \in[0,1]$ such that $ (\hat{S}^R )^{-1}(a) \succ
^{U,a,b} t_{ab}$ and $ (\hat{S}^R )^{-1}(b) \succ^{U,b,a}
t_{ba}$, implying that $\hat{S}^R(i) = \hat{\sigma}(\xi_i)$ for some
$\hat{\sigma} \in\mathcal{Q}_U$. As parallel arguments hold for
$\hat
{T}^Q$, the statement of the lemma follows.
\end{pf}

\subsection{\texorpdfstring{Proof of Lemma \protect\ref{le:rademacher1.easy}}
{Proof of Lemma 5.6}}
\label{le:rademacher1.easypf}

Recall the definition of $\delta_{UV}$ from Section~\ref{sec:K2Case},
which we can now interpret for arbitrary $K$ according to the
definitions of $\mathcal{Q}_U$ and $\mathcal{Q}_V$ in Section~\ref{sec:QU and QV} above. We use a symmetrization argument of Hoeffding
[\citet{hoeffding1963probability,clemenccon2008ranking}] to bound
$\mathbb{E} (\max_{(Q,R) \in\mathcal{Q}_\mu^m
\times\mathcal{Q}_\nu^n} \delta_{UV} )$. Let $\mathcal{M}_\mathcal
{I}$ denote the set of permutations of $1,\ldots,m$ which map
$1,\ldots
,m-|\mathcal{I}|$ to $i \notin\mathcal{I}$, and let $\mathcal
{M}_\mathcal{J}$ be defined analogously for permutations on $1,\ldots
,n$. Let $\mathcal{M} = \mathcal{M}_\mathcal{I} \times\mathcal
{M}_\mathcal{J}$ and let $Z = |\mathcal{M}|$. Let $\xi',\zeta'$ be
identically distributed as $\xi$ and $\zeta$, and independent of $U$
and $V$. Let $\xi_\mathcal{I}$ and $\zeta_\mathcal{J}$ be defined
as in
Section~\ref{sec:QU and QV}. To abbreviate the notation, let
$g_{\sigma
\tau}(x,y) = \omega(x,y)\Gamma_{\sigma(x)\tau(y)}$, and let
$\mathcal
{Q} = \mathcal{Q}_\mu^m \times\mathcal{Q}_\nu^n \times\mathcal{Q}_U
\times\mathcal{Q}_V$. It holds for $(Q,R) \in\mathcal{Q}_\mu^m
\times
\mathcal{Q}_\nu^n$ that
\[
\mathbb{E} \Bigl( \max_{Q,R} \delta_{UV} \Bigr) =
\mathbb{E} \Bigl( \sup_{(Q,R,\sigma,\tau)\in
\mathcal
{Q}} \bigl|
G_{\sigma,\tau}(\xi,\zeta) - \mathbb{E} \bigl(G_{\sigma\tau}
\bigl(\xi',\zeta' \bigr)\bigr|U,V \bigr) \bigr| \vert
\xi_\mathcal{I},\zeta_\mathcal{J} \Bigr),
\]
which by convexity can be upper-bounded by
\begin{eqnarray*}
&& \mathbb{E} \Bigl(\sup_{(Q,R,\sigma,\tau)\in
\mathcal
{Q}}\bigl |
G_{\sigma,\tau}(\xi,\zeta) - G_{\sigma\tau} \bigl(\xi',
\zeta' \bigr) \bigr| \vert\xi_\mathcal{I},\zeta_\mathcal{J}
\Bigr)
\\
&&\qquad = \mathbb{E} \biggl( \sup_{(Q,R,\sigma,\tau)\in
\mathcal{Q}} \biggl
\llvert\frac{1}{mn} \sum_{i\notin\mathcal{I}} \sum
_{j
\notin\mathcal{J}} g_{\sigma\tau}(\xi_i,
\zeta_j)  - g_{\sigma\tau
} \bigl(\xi_i',
\zeta_j' \bigr) \biggr\rrvert\Big\vert
\xi_\mathcal{I}, \zeta_\mathcal{J} \biggr)
\\
& &\qquad= \mathbb{E} \Biggl( \sup_{(Q,R,\sigma,\tau)\in
\mathcal{Q}} \Biggl
\llvert\frac{|\overline{\mathcal{I}}| |\overline{\mathcal
{J}}|}{Zmn} \sum_{\pi,\eta\in\mathcal{M}}
\frac{1}{\ell} \sum_{i=1}^\ell
g_{\sigma\tau}(\xi_{\pi(i)},\zeta_{\eta(j)}) \\
&&\hspace*{185pt}{}- g_{\sigma
\tau}
\bigl(\xi_{\pi(i)}',\zeta_{\eta(j)}'
\bigr) \Biggr\rrvert\Big\vert\xi_\mathcal{I},\zeta_\mathcal{J} \Biggr),
\end{eqnarray*}
{since the permutations $\pi$ and $\eta$ weight each $(i,j)$ term
equally for $i \notin\mathcal{I}$ and $j\notin\mathcal{J}$; by
convexity again, and then linearity of expectation, we have}
\begin{eqnarray*}
& \leq&\mathbb{E}\Biggl( \frac{|\overline{\mathcal{I}}|
|\overline{\mathcal{J}}|}{Zmn} \sum
_{\pi,\eta\in\mathcal{M}} \sup_{(Q,R,\sigma,\tau)\in\mathcal
{Q}} \Biggl\llvert
\frac{1}{\ell} \sum_{i=1}^\ell
g_{\sigma\tau}(\xi_{\pi(i)},\zeta_{\eta(j)}) - g_{\sigma
\tau}
\bigl(\xi_{\pi(i)}',\zeta_{\eta(j)}'
\bigr) \Biggr\rrvert\Big\vert\xi_\mathcal{I},\zeta_\mathcal{J} \Biggr)
\\
& =&\frac{|\overline{\mathcal{I}}| |\overline{\mathcal{J}}|}{mn} \mathbb{E} \Biggl(\sup_{(Q,R,\sigma,\tau)\in
\mathcal{Q}}
\Biggl\llvert\frac{1}{\ell} \sum_{i=1}^\ell
g_{\sigma\tau}(\xi_{i},\zeta_{i})- g_{\sigma\tau}
\bigl(\xi_{i}',\zeta_{i}'
\bigr) \Biggr\rrvert\Big\vert\xi_\mathcal{I},\zeta_\mathcal{J} \Biggr).
\end{eqnarray*}

We may now introduce independent and identically distributed Rademacher
variables $r_1,\ldots,r_\ell$, and use standard Rademacher
symmetrization arguments [see, e.g., \citet{bousquet2004introduction}]
to show that the final quantity above is equal to
\begin{eqnarray*}
& &\frac{|\overline{\mathcal{I}}| |\overline{\mathcal{J}}|}{mn} \mathbb{E} \Biggl(\sup_{(Q,R,\sigma,\tau)\in
\mathcal{Q}}
\Biggl\llvert\frac{1}{\ell} \sum_{i=1}^\ell
r_i \bigl(g_{\sigma\tau}(\xi_{i},
\zeta_{i})- g_{\sigma\tau} \bigl(\xi_{i}',
\zeta_{i}' \bigr) \bigr) \Biggr\rrvert\Big\vert
\xi_\mathcal{I},\zeta_\mathcal{J} \Biggr)
\\
&&\qquad \leq\frac{|\overline{\mathcal{I}}| |\overline{\mathcal
{J}}|}{mn} \mathbb{E} \Biggl(\sup
_{(Q,R,\sigma,\tau)\in
\mathcal{Q}} \Biggl\llvert\frac{1}{\ell} \sum
_{i=1}^\ell r_i g_{\sigma\tau
}(
\xi_{i},\zeta_{i}) \Biggr\rrvert+ \Biggl\llvert
\frac{1}{\ell} \sum_{i=1}^\ell
r_i g_{\sigma\tau} \bigl(\xi_{i}',
\zeta_{i}' \bigr) \Biggr\rrvert\Big\vert\xi_\mathcal
{I},\zeta_\mathcal{J} \Biggr)
\\
& &\qquad\leq2 \frac{|\overline{\mathcal{I}}| |\overline{\mathcal
{J}}|}{mn} \mathbb{E} \Biggl(\sup
_{(Q,R,\sigma,\tau)\in
\mathcal{Q}} \Biggl\llvert\frac{1}{\ell} \sum
_{i=1}^\ell r_i g_{\sigma\tau
}(
\xi_{i},\zeta_{i}) \Biggr\rrvert\Big\vert
\xi_\mathcal{I},\zeta_\mathcal{J} \Biggr).
\end{eqnarray*}

To bound this expectation, note that for fixed $\mathcal{I}, \mathcal
{J},Q,R$ (inducing a fixed $U$ and $V$), and fixed $(\sigma,\tau) \in
\mathcal{Q}_U \times\mathcal{Q}_V$, a Hoeffding inequality gives
%
\begin{equation}
\label{eq:hoeffding} \mathbb{P} \Biggl( \Biggl\llvert\frac
{1}{\ell} \sum
_{i=1}^\ell r_i g_{\sigma\tau
}(
\xi_i,\zeta_j) \Biggr\rrvert\geq\varepsilon\Big\vert
\xi_\mathcal{I},\zeta_\mathcal{J} \Biggr) \leq2e^{-2\ell\varepsilon^2}.
\end{equation}

We may now apply \eqref{eq:hoeffding} in conjunction with a union bound
over all $(Q,R,\sigma,\tau)\in\mathcal{Q}$ as follows. For fixed
$Q,R,a,b$, the set $\{i\dvtx \xi_i \succeq^{U,a,b} t_{ab}\}$ can be chosen
at most $\ell+1$ ways by varying $t_{ab}$. As a result, the set $\xi
_1,\ldots,\xi_\ell$ can be partitioned at most $(\ell+1)^{K \choose2}$
ways by varying $\sigma\in\mathcal{Q}_U$. Analogously, the set
$\zeta
_1,\ldots,\zeta_\ell$ can be partitioned the same number of ways by
varying $\tau\in\mathcal{Q}_V$. For fixed $\mathcal{I}, \mathcal{J}$,
the functions $U$ and $V$ can be chosen $K^{|\mathcal{I}|+|\mathcal
{J}|}$ different ways by varying $Q$ and $R$. Hence, a union bound gives
\[
\mathbb{P} \Biggl( \sup_{(Q,R,\sigma,\tau)\in\mathcal{Q}} \Biggl\llvert
\frac
{1}{\ell}\sum_{i=1}^\ell
r_i g_{\sigma\tau}(\xi_i \zeta_i)
\Biggr\rrvert\geq\varepsilon\Big\vert\xi_\mathcal{I},\zeta_\mathcal{J}
\Biggr) \leq K^{|\mathcal{I}| + |\mathcal{J}|} (\ell+ 1)^{2{K
\choose2}} \cdot2e^{-2\ell\varepsilon^2}.
\]
Since this expression is of the form $\mathbb{P}(X \geq t) \leq f(t)$
for $X$ nonnegative, we may apply the inequality $\mathbb
{E}(X) \leq\int_0^\infty\min\{ 1, f(t) \} \,dt$ to yield
\begin{eqnarray*}
&&2 \frac{|\overline{\mathcal{I}}| |\overline{\mathcal{J}}|}{mn} \mathbb{E} \Biggl(\sup_{(Q,R,\sigma, \tau) \in
\mathcal
{Q}}
\Biggl\llvert\frac{1}{\ell} \sum_{i=1}^\ell
r_i g_{\sigma\tau}(\xi_i,\zeta_i)
\Biggr\rrvert\Big\vert\xi_\mathcal{I},\zeta_\mathcal{J} \Biggr)
\\
&&\qquad\leq4\sqrt{\frac{(|\mathcal{I}|+|\mathcal{J}|) \log K + 2{K
\choose
2}\log(\ell+1) + \log2}{2\ell}}.
\end{eqnarray*}
Since the bound holds for any $\xi_\mathcal{I}, \zeta_\mathcal{J}$, the
same bound holds when the conditioning is removed and $\xi_\mathcal
{I},\zeta_\mathcal{J}$ are chosen randomly, thus proving the lemma.

\subsection{\texorpdfstring{Proof of Lemma \protect\ref{le:rademacher2.easy}}
{Proof of Lemma 5.7}}
\label{le:rademacher2.easypf}

To abbreviate notation, let $\mathcal{Q} = \mathcal{Q}_\nu^n \times
\mathcal{Q}_U$. Let $r_1,\ldots,r_m$ be Rademacher variables as in the
proof of Lemma~\ref{le:rademacher1.easy}. By a standard Rademacher
symmetrization,
\begin{eqnarray*}
&&\mathbb{E} \Biggl(\sup_{(R,\sigma) \in\mathcal
{Q}} \Biggl\{ \max
_{1\leq a \leq K} \Biggl\llvert\sigma^{-1}(a) -
\frac{1}{m} \sum_{i=1}^m 1 \bigl\{
\sigma(\xi_i)=a \bigr\} \Biggr\rrvert \Biggr\} \Big\vert
\zeta_\mathcal{J} \Biggr)
\\
&&\qquad\leq2\mathbb{E} \Biggl( \sup_{(R,\sigma) \in
\mathcal
{Q}} \Biggl\{
\max_{1\leq a \leq K} \Biggl\llvert\frac{1}{m} \sum
_{i=1}^m r_i 1 \bigl\{\sigma(
\xi_i)=a \bigr\} \Biggr\rrvert \Biggr\} \Big\vert\zeta_\mathcal{J}
\Biggr).
\end{eqnarray*}
As in the proof of Lemma~\ref{le:rademacher1.easy}, a Hoeffding
inequality and union bound yield
\begin{eqnarray*}
&&\mathbb{P} \Biggl( \sup_{(R,\sigma) \in\mathcal{Q}} \Biggl\{ \max
_{1\leq
a \leq K} \Biggl\llvert\sigma^{-1}(a) -
\frac{1}{m} \sum_{i=1}^m 1 \bigl\{
\sigma(\xi_i) = a \bigr\} \Biggr\rrvert \Biggr\} \geq\varepsilon\Big\vert
\zeta_\mathcal{J} \Biggr)\\
&&\qquad \leq K^{|\mathcal{J}|} (m+1)^{K\choose2} K
\cdot2e^{-2m\varepsilon^2},
\end{eqnarray*}
and applying $\mathbb{E}(|X|) \leq\int_0^\infty\min
\{ 1, f(t) \} \,dt$ for $\mathbb{P}(|X| \geq t) \leq f(t)$
then gives
\begin{eqnarray*}
&&2 \mathbb{E} \Biggl(\sup_{(R,\sigma) \in\mathcal
{Q}} \Biggl\{ \max
_{1\leq a \leq K} \Biggl\llvert\sigma^{-1}(a) -
\frac{1}{m} \sum_{i=1}^m 1 \bigl\{
\sigma(\xi_i)=a \bigr\} \Biggr\rrvert \Biggr\} \Big\vert
\zeta_\mathcal{J} \Biggr)
\\
&&\qquad\leq4 \sqrt{\frac{(|\mathcal{J}|+1)\log K + {K\choose2} \log(m+1)
+ \log2}{2m}}.
\end{eqnarray*}
As in the proof of Lemma~\ref{le:rademacher1.easy}, removing the
conditioning on $\zeta_\mathcal{J}$ does not alter the bound. Parallel
arguments apply to $\tau\in\mathcal{Q}_V$, and the lemma follows.

\subsection{\texorpdfstring{Covering argument to establish Theorem \protect\ref{th:main}}
{Covering argument to establish Theorem 4.1}}
\label{sec:epsCovering}

The establishment of Theorem~\ref{th:main} from Proposition~\ref
{prop:hyperplane} proceeds as follows. For $\mathcal{F} \subset
[0,1]^{K \times K}$, recall that $h_{\mathcal{F}}(\Gamma) = \sup_{F
\in
\mathcal{F}} \langle\Gamma,F \rangle = \sup_{F \in\mathcal{F}}
\operatorname
{tr}( \Gamma^T F ) $. By the Cauchy--Schwarz inequality, $h_\mathcal
{F}$ is Lipschitz continuous,
\[
\bigl\llvert h_{\mathcal{F}}(\Gamma) - h_{\mathcal{F}} \bigl(
\Gamma' \bigr) \bigr\rrvert\leq\sup_{F \in\mathcal{F}} \bigl
\llvert\bigl\langle\Gamma-\Gamma',F \bigr\rangle \bigr\rrvert\leq K
\bigl\| \Gamma- \Gamma'\bigr\|.
\]
Let $\mathcal{B}_\varepsilon$ denote an $\varepsilon$-cover in $\| \cdot\|$
for $[-1,1]^{K \times K}$, with $\Gamma^\mathcal{B}$ the closest point
in $\mathcal{B}_\varepsilon$ to a given $\Gamma$. The triangle inequality,
Lipschitz condition and $\mathcal{B}_\varepsilon$ imply
\begin{eqnarray*}
&&\sup_{\Gamma\in[-1,1]^{K \times K}}  \bigl\llvert h_{{\mathcal
{F}_{\mu\nu}^A}}(\Gamma) -
h_{{\mathcal{F}_{\mu\nu}^\omega}}(\Gamma) \bigr\rrvert\\
&&\qquad\leq\sup_{\Gamma\in[-1,1]^{K \times K}} \bigl\{
\bigl\llvert h_{{\mathcal{F}_{\mu\nu}^A}
}(\Gamma) - h_{{\mathcal{F}_{\mu\nu}^A}} \bigl(
\Gamma^\mathcal{B} \bigr) \bigr\rrvert
\\
&&\hspace*{55pt}\qquad\quad{} + \bigl\llvert h_{{\mathcal{F}_{\mu\nu}^A}} \bigl(\Gamma ^\mathcal{B} \bigr) -
h_{{\mathcal{F}_{\mu\nu}^\omega}
} \bigl(\Gamma^\mathcal{B} \bigr) \bigr\rrvert+ \bigl
\llvert h_{{\mathcal{F}_{\mu\nu}^\omega}} \bigl(\Gamma^\mathcal{B} \bigr) -
h_{{\mathcal{F}_{\mu\nu}^\omega}}(\Gamma) \bigr\rrvert \bigr\}
\\
&&\qquad \leq\sup_{\Gamma\in[-1,1]^{K \times K}} \bigl\{ \bigl\llvert h_{{\mathcal{F}_{\mu\nu}^A}
}
\bigl( \Gamma^\mathcal{B} \bigr) - h_{{\mathcal{F}_{\mu\nu}^\omega
}} \bigl(
\Gamma^\mathcal{B} \bigr) \bigr\rrvert+ 2 K \bigl\| \Gamma- \Gamma^\mathcal{B}
\bigr\| \bigr\}
\\
&&\qquad \leq\sup_{\Gamma\in[-1,1]^{K \times K}} \bigl\llvert h_{{\mathcal{F}_{\mu\nu}^A}} \bigl(\Gamma
^\mathcal{B} \bigr) - h_{{\mathcal{F}_{\mu\nu}^\omega}} \bigl(\Gamma^\mathcal{B}
\bigr) \bigr\rrvert+ 2 K \varepsilon
\\
&&\qquad = \max_{\Gamma\in\mathcal{B}_\varepsilon} \bigl\llvert h_{{\mathcal{F}_{\mu\nu}^A}}(\Gamma) -
h_{{\mathcal{F}_{\mu\nu}^\omega}}(\Gamma) \bigr\rrvert+ 2 K \varepsilon.
\end{eqnarray*}

Now let $C_K$ and $n_K$ be defined as in Proposition~\ref
{prop:hyperplane}, and set $\varepsilon= C_K / n^{1/4}$. It follows by
the above relation, a union bound, and Proposition~\ref
{prop:hyperplane} that
\begin{eqnarray*}
&&  \mathbb{P} \Bigl(\max_{(\mu,\nu) \in\Omega_{\rho
n} \times
\Omega_n} \Bigl\{ \sup
_{\Gamma\in[-1,1]^{K \times K}} \bigl| h_{{\mathcal{F}_{\mu\nu}^A}
}(\Gamma) - h_{{\mathcal{F}_{\mu\nu}^\omega}}(\Gamma)
\bigr| \Bigr\} \geq3\varepsilon \Bigr)
\\
&&\qquad \leq\mathbb{P} \Bigl(\max_{(\mu,\nu) \in\Omega_{\rho n}
\times
\Omega_n} \Bigl\{ \max
_{\Gamma\in\mathcal{B}_{\varepsilon/K}} \bigl| h_{{\mathcal{F}_{\mu\nu
}^A}}(\Gamma) - h_{{\mathcal{F}_{\mu\nu}^\omega}}(\Gamma)
\bigr| \Bigr\} \geq\varepsilon \Bigr)
\\
&&\qquad \leq\sum_{(\mu,\nu) \in\Omega_{\rho n} \times\Omega_n} \sum_{\Gamma\in\mathcal{B}_{\varepsilon/K}}
\mathbb{P} \bigl( \bigl| h_{{\mathcal{F}_{\mu\nu}^A}
}(\Gamma) - h_{{\mathcal{F}_{\mu\nu}^\omega}}(\Gamma) \bigr| \geq
\varepsilon \bigr)
\\
&&\qquad \leq|\Omega_{\rho n}| |\Omega_n| |\mathcal{B}_{\varepsilon/K}|
2e^{-\sqrt{n} [2\rho/(\rho+1)]} \bigl[1+o(1) \bigr]
\end{eqnarray*}
for all $n \geq n_K$. The result of Theorem~\ref{th:main} then follows,
since we have that $|\Omega_n| =
{n+K-1 \choose K-1}$, and $\mathcal{B}_{\varepsilon/K}$ can be chosen such
that $|\mathcal{B}_{\varepsilon/K}| \leq(1 + K^2 / \varepsilon)^{K^2}$.

\section{\texorpdfstring{Statement and Proof of Lemma
\lowercase{\protect\ref{le:sim1}}}
{Statement and Proof of Lemma C.1}}

\label{sec:simulationLemmaPf}

To evaluate the excess risk quantities reported in Section~\ref{subsec:simulations}, we require both that $\phi^* =
\operatorname
{argmax}_{\phi\in\Phi} L_{\rho_n \omega_\beta}(\phi)$ exist, and that
$L_{\rho_n \omega_\beta}(\phi)$ be computable. The following lemma
establishes this, using the fact that each $\rho_n \omega_\beta
(x,y)$ is a separable function plus a constant. Given a triple $\phi=
(\mu,\nu,\theta)$ of two-class blockmodel parameters, it shows that
$L_{\rho_n\omega_\beta}(\phi)$, which nominally involves an
optimization over all measure-preserving maps of $[0,1]$, can be
reduced to a maximization over four cases, and thus evaluated tractably.

\begin{lemma}
\label{le:sim1}
Given $\mu\in[0,1]^2$, let $\sigma^{(1)}, \sigma^{(2)} \in\mathcal
{Q}_\mu$ denote the mappings
\[
\sigma^{(1)}(x) = \cases{%
 1, &\quad
$\mbox{if }0 \leq x < \mu_1,$
\vspace*{2pt}\cr
2, &\quad $\mbox{if }\mu_1 \leq x \leq1;$ }\qquad
\sigma^{(2)}(x) = \cases{ %
1, &\quad
$\mbox{if }1 - \mu_2 \leq x \leq1,$\vspace*{2pt}\cr
2, & \quad $\mbox{if }0 \leq x < 1 - \mu_2;$}
\]
and let $\tau^{(1)}, \tau^{(2)} \in
\mathcal{Q}_\nu$ be defined analogously, given $\nu
\in[0,1]^2$. Given $\phi= (\mu,\nu,\theta)$, terms $L_{\rho_n
\omega_\beta}(
\phi)$ and $R_{\rho_n \omega_\beta}(\phi)$ from Theorem~\ref
{th:learning} equal$ $
\begin{eqnarray*}
L_{\rho_n \omega_\beta}(\phi) & =& \max_{(i,j) \in\{1,2\}^2} \sum
_{a=1}^2 \sum_{b=1}^2
\bigl(\rho_n \omega/\sigma^{(i)} \tau^{(j)}
\bigr)_{ab} \log \biggl(\frac{\theta_{ab}}{1-\theta_{ab}} \biggr) \\
&&\hspace*{75pt}{}+ \mu_a
\nu_b \log(1-\theta_{ab}),
\\
R_{\rho_n \omega_\beta}(\phi) & =& \min_{(i,j) \in\{1,2\}^2} \Biggl\{ \sum
_{a=1}^2 \sum_{b=1}^2
-2 \bigl( \rho_n \omega/ \sigma^{(i)} \tau^{(j)}
\bigr)_{ab} \theta_{ab} + \mu_a
\nu_b \theta_{ab}^2 \Biggr\} \\
&&{}+ \int
_{[0,1]^2} \omega(x,y)^2 \,dx \,dy.
\end{eqnarray*}
\end{lemma}
\begin{pf} Below we establish the claimed expression for $L_{\rho_n
\omega_\beta
}(
\phi)$; analogous arguments yield the result for $R_{\rho_n \omega
_\beta}(\phi)$. First,
define $L_\omega(\phi; \sigma,\tau)$ as
\[
L_\omega( \phi;\sigma,\tau) = \sum_{a=1}^2
\sum_{b=1}^2 (\omega/\sigma
\tau)_{ab} \log \biggl(\frac{\theta_{ab}}{1-\theta_{ab}} \biggr) + \mu_a
\nu_b \log(1-\theta_{ab}).
\]
Next, let $\sigma^* \vert \tau= \operatorname{argmax}_{\sigma\in
\mathcal{Q}_\mu} L_\omega(\phi;\sigma,\tau)$, with the convention
that\break
$\operatorname{argmax}_{\sigma\in\mathcal{Q}_\mu} (\cdot)$ is
undefined if no maximizer exists. We then see that
\begin{eqnarray}
\sigma^* \vert\tau& =& \mathop{\operatorname{argmax}}_{\sigma\in\mathcal
{Q}_\mu} \sum
_{a=1}^2 \sum_{b=1}^2
(\omega/\sigma\tau)_{ab} \log \biggl( \frac{\theta_{ab}}{1-\theta
_{ab}} \biggr)
\nonumber\\[-2pt]
& =& \mathop{\operatorname{argmax}}_{\sigma\in\mathcal{Q}_\mu} \sum_{a=1}^2
\sum_{b=1}^2 \int_{\sigma^{-1}(a) \times\tau^{-1}(b)}
\omega(x,y) \,dx \,dy \log \biggl( \frac{\theta
_{ab}}{1-\theta_{ab}} \biggr)
\nonumber\\[-2pt]
& =& \mathop{\operatorname{argmax}}_{\sigma\in\mathcal{Q}_\mu} \sum_{a=1}^2
\int_{\sigma^{-1}(a)} g_a(x) \,dx\nonumber\\[-2pt]
\eqntext{\mbox{with }
\displaystyle g_a(x)
= \sum_{b=1}^2 \int
_{\tau^{-1}(b)} \omega(x,y) \log \biggl(\frac{\theta
_{ab}}{1-\theta_{ab}} \biggr) \,dy}
\\[-2pt]
& =& \mathop{\operatorname{argmax}}_{\sigma\in\mathcal{Q}_\mu} \int_{[0,1]}
g_2(x) \,dx + \int_{[0,1]} \bigl(g_1(x)
- g_2(x) \bigr) 1 \bigl\{\sigma(x)=1 \bigr\} \,dx.\nonumber
\end{eqnarray}
It can be seen that $\sigma^* \vert \tau$ is always defined and
assigns the $\mu_1$-quantile of $g_1(x) - g_2(x)$ to class $1$. Since
$\rho_n \omega_\beta(x,y) = \rho_n ( f_\beta(x)f_\beta(y) + 1/2)$,
$g_1(x) - g_2(x)$ is affine in $f_\beta(x)$, and can be written as $m
f_\beta(x)+c$ for some scalars $m$ and $c$. As $f_\beta$ is monotone,
the $\mu_1$-quantile will either be $[0,\mu_1]$ or $[1-\mu
_1,1]$---depending on the sign of $m$---meaning that $\sigma^* \vert
\tau$ equals either $\sigma^{(1)}$ or $\sigma^{(2)}$ for any $\tau$.
Analogously, $\tau^* \vert \sigma$ equals either $\tau^{(1)}$ or
$\tau^{(2)}$ for any $\sigma$. Hence,
\begin{eqnarray*}
L_{\rho_n \omega_\beta}(\phi) &=& \sup_{\sigma\in\mathcal{Q}_\mu
} \sup
_{\tau\in\mathcal{Q}_\nu} L_{\rho_n \omega_\beta}(\phi;\sigma,\tau)
\\[-2pt]
& =& \sup_{\sigma\in\mathcal{Q}_\mu} L_{\rho_n \omega_\beta} \bigl(\phi; \sigma, \bigl(
\tau^* \vert\sigma \bigr) \bigr) \leq\sup_{\sigma\in\mathcal
{Q}_\mu}
L_{\rho_n \omega_\beta} \bigl(\phi; \bigl(\sigma^* \vert \bigl(\tau^* \vert\sigma
\bigr) \bigr), \bigl(\tau^* \vert\sigma \bigr) \bigr)
\\[-2pt]
&  =& \max_{(i,j) \in\{1,2\}^2} L_{\rho_n \omega_\beta
} \bigl(\phi;
\sigma^{(i)}, \tau^{(j)} \bigr).
\end{eqnarray*}
Thus $L_{\rho_n \omega_\beta}(\phi)$, which nominally involves a
supremum over every pair $(\sigma, \tau) \in\mathcal{Q}_\mu\times
\mathcal{Q}_\nu$, is reduced to a maximization over $\sigma^{\smash
{(1)}},\sigma^{\smash{(2)}}$ and $\tau^{\smash{(1)}},\tau^{\smash{(2)}}$.
\end{pf}

\begin{corollary}\label{cor:sim1}
The quantity $\sup_{\phi\in{\Phi}} L_{\rho_n \omega_\beta
}(\phi
)$ is achieved by $\phi^{(ij)} = (\mu,\nu,\omega/\sigma^{(i)}\tau
^{(j)})$ for some $(\mu,\nu) \in\Omega_m \times\Omega_n$ and $(i,j)
\in\{1,2\}^2$.
\end{corollary}
\begin{pf}
For any $\phi= (\mu,\nu,\theta) \in{\Phi}$, it holds that
\begin{eqnarray*}
L_{\rho_n \omega_\beta} (\phi) & =& \max_{(i,j) \in\{1,2\}^2} L_{\rho_n
\omega_\beta}
\bigl(\phi;\sigma^{(i)},\tau^{(j)} \bigr)
\\[-2pt]
& \leq&\max_{(i,j) \in\{1,2\}^2} L_{\rho_n \omega_\beta} \bigl(\phi^{(ij)};
\sigma^{(i)}, \tau^{(j)} \bigr)
\\[-2pt]
& = &\max_{(i,j) \in\{1,2\}^2} L_{\rho_n \omega_\beta} \bigl(\phi^{(ij)}
\bigr),
\end{eqnarray*}
where the first line holds by Lemma~\ref{le:sim1}, the second because
$p \log x + (1-p)\log(1-x)$ is maximized over $0\leq x\leq1$ by $x=p$,
and the third by the definition of $L_{\rho_n \omega_\beta}(\cdot;\cdot,\cdot)$.
\end{pf}
\end{appendix}
\section*{Acknowledgment}

The first author wishes to thank Peter Bickel for helpful advice and
feedback.

%



\printaddresses


\begin{thebibliography}{29}

\bibitem[\protect\citeauthoryear{Airoldi et~al.}{2008}]{airoldi2008mixed}
%
\begin{barticle}[author]
\bauthor{\bsnm{Airoldi},~\bfnm{E.~M.}\binits{E.~M.}},
\bauthor{\bsnm{Blei},~\bfnm{D.~M.}\binits{D.~M.}},
\bauthor{\bsnm{Fienberg},~\bfnm{S.~E.}\binits{S.~E.}} \AND
\bauthor{\bsnm{Xing},~\bfnm{E.~P.}\binits{E.~P.}}
(\byear{2008}).
\btitle{Mixed membership stochastic blockmodels}.
\bjournal{J. Mach. Learn. Res.}
\bvolume{9}
\bpages{1981--2014}.
\bptok{imsref}%
\end{barticle}
%
\endbibitem

\bibitem[\protect\citeauthoryear{Alon et~al.}{2003}]{alon2003random}
%
\begin{barticle}[mr]
\bauthor{\bsnm{Alon},~\bfnm{Noga}\binits{N.}}, \bauthor{\bparticle
{Fernandez
de~la} \bsnm{Vega},~\bfnm{W.}\binits{W.}},
\bauthor{\bsnm{Kannan},~\bfnm{Ravi}\binits{R.}} \AND
\bauthor{\bsnm{Karpinski},~\bfnm{Marek}\binits{M.}}
(\byear{2003}).
\btitle{Random sampling and approximation of {MAX}-{CSP}s}.
\bjournal{J. Comput. System Sci.}
\bvolume{67}
\bpages{212--243}.
\bid{doi={10.1016/S0022-0000(03)00008-4}, issn={0022-0000}, mr={2022830}}
\bptok{imsref}%
\end{barticle}
%
\endbibitem

\bibitem[\protect\citeauthoryear{Bickel and
Chen}{2009}]{bickel2009nonparametric}
%
\begin{barticle}[author]
\bauthor{\bsnm{Bickel},~\bfnm{P.~J.}\binits{P.~J.}} \AND
\bauthor{\bsnm{Chen},~\bfnm{A.}\binits{A.}}
(\byear{2009}).
\btitle{A nonparametric view of network models and Newman--Girvan and other
modularities}.
\bjournal{Proc. Natl. Acad. Sci. USA}
\bvolume{106}
\bpages{21068--21073}.
\bptok{imsref}%
\end{barticle}
%
\endbibitem

\bibitem[\protect\citeauthoryear{Bickel, Chen and
Levina}{2011}]{bickelmethod}
%
\begin{barticle}[mr]
\bauthor{\bsnm{Bickel},~\bfnm{Peter~J.}\binits{P.~J.}},
\bauthor{\bsnm{Chen},~\bfnm{Aiyou}\binits{A.}} \AND
\bauthor{\bsnm{Levina},~\bfnm{Elizaveta}\binits{E.}}
(\byear{2011}).
\btitle{The method of moments and degree distributions for network models}.
\bjournal{Ann. Statist.}
\bvolume{39}
\bpages{2280--2301}.
\bid{doi={10.1214/11-AOS904}, issn={0090-5364}, mr={2906868}}
\bptok{imsref}%
\end{barticle}
%
\endbibitem

\bibitem[\protect\citeauthoryear{Borgs et~al.}{2006}]{borgs2006graph}
%
\begin{bincollection}[mr]
\bauthor{\bsnm{Borgs},~\bfnm{Christian}\binits{C.}},
\bauthor{\bsnm{Chayes},~\bfnm{Jennifer}\binits{J.}},
\bauthor{\bsnm{Lov{\'a}sz},~\bfnm{L{\'a}szl{\'o}}\binits{L.}},
\bauthor{\bsnm{S{\'o}s},~\bfnm{Vera~T.}\binits{V.~T.}},
\bauthor{\bsnm{Szegedy},~\bfnm{Bal{\'a}zs}\binits{B.}} \AND
\bauthor{\bsnm{Vesztergombi},~\bfnm{Katalin}\binits{K.}}
(\byear{2006}).
\btitle{Graph limits and parameter testing}.
In \bbooktitle{S{TOC}'06: {P}roceedings of the 38th {A}nnual {ACM} {S}ymposium
on {T}heory of {C}omputing}
\bpages{261--270}.
\bpublisher{ACM}, \blocation{New York}.
\bid{doi={10.1145/1132516.1132556}, mr={2277152}}
\bptok{imsref}%
\end{bincollection}
%
\endbibitem

\bibitem[\protect\citeauthoryear{Borgs et~al.}{2008}]{borgs2008convergent}
%
\begin{barticle}[mr]
\bauthor{\bsnm{Borgs},~\bfnm{C.}\binits{C.}},
\bauthor{\bsnm{Chayes},~\bfnm{J.~T.}\binits{J.~T.}},
\bauthor{\bsnm{Lov{\'a}sz},~\bfnm{L.}\binits{L.}},
\bauthor{\bsnm{S{\'o}s},~\bfnm{V.~T.}\binits{V.~T.}} \AND
\bauthor{\bsnm{Vesztergombi},~\bfnm{K.}\binits{K.}}
(\byear{2008}).
\btitle{Convergent sequences of dense graphs. {I}. {S}ubgraph frequencies,
metric properties and testing}.
\bjournal{Adv. Math.}
\bvolume{219}
\bpages{1801--1851}.
\bid{doi={10.1016/j.aim.2008.07.008}, issn={0001-8708}, mr={2455626}}
\bptok{imsref}%
\end{barticle}
%
\endbibitem

\bibitem[\protect\citeauthoryear{Borgs et~al.}{2012}]{borgs2007convergent}
%
\begin{barticle}[mr]
\bauthor{\bsnm{Borgs},~\bfnm{C.}\binits{C.}},
\bauthor{\bsnm{Chayes},~\bfnm{J.~T.}\binits{J.~T.}},
\bauthor{\bsnm{Lov{\'a}sz},~\bfnm{L.}\binits{L.}},
\bauthor{\bsnm{S{\'o}s},~\bfnm{V.~T.}\binits{V.~T.}} \AND
\bauthor{\bsnm{Vesztergombi},~\bfnm{K.}\binits{K.}}
(\byear{2012}).
\btitle{Convergent sequences of dense graphs. {II}. {M}ultiway cuts and
statistical physics}.
\bjournal{Ann. of Math. (2)}
\bvolume{176}
\bpages{151--219}.
\bid{doi={10.4007/annals.2012.176.1.2}, issn={0003-486X}, mr={2925382}}
\bptok{imsref}%
\end{barticle}
%
\endbibitem

\bibitem[\protect\citeauthoryear{Bousquet, Boucheron and
Lugosi}{2004}]{bousquet2004introduction}
%
\begin{bincollection}[author]
\bauthor{\bsnm{Bousquet},~\bfnm{O.}\binits{O.}},
\bauthor{\bsnm{Boucheron},~\bfnm{S.}\binits{S.}} \AND
\bauthor{\bsnm{Lugosi},~\bfnm{G.}\binits{G.}}
(\byear{2004}).
\btitle{Introduction to statistical learning theory}.
In \bbooktitle{Advanced Lectures on Machine Learning}
(\beditor{\bfnm{O.}\binits{O.}~\bsnm{Bousquet}},
\beditor{\bfnm{U.}\binits{U.}~\bparticle{von} \bsnm{Luxburg}} \AND
\beditor{\bfnm{G.}\binits{G.}~\bsnm{R{\"a}tsch}}, eds.)
\bpages{169--207}.
\bpublisher{Springer}, \blocation{Berlin}.
\bptok{imsref}%
\end{bincollection}
%
\endbibitem

\bibitem[\protect\citeauthoryear{Chatterjee}{2012}]{chatterjee2012matrix}
%
\begin{bmisc}[author]
\bauthor{\bsnm{Chatterjee},~\bfnm{S.}\binits{S.}}
(\byear{2012}).
\bhowpublished{Matrix estimation by universal singular value thresholding.
Preprint. Available at \arxivurl{arXiv:1212.1247}.}
\bptok{imsref}%
\end{bmisc}
%
\endbibitem

\bibitem[\protect\citeauthoryear{Choi, Wolfe and
Airoldi}{2012}]{choi2012stochastic}
%
\begin{barticle}[mr]
\bauthor{\bsnm{Choi},~\bfnm{D.~S.}\binits{D.~S.}},
\bauthor{\bsnm{Wolfe},~\bfnm{P.~J.}\binits{P.~J.}} \AND
\bauthor{\bsnm{Airoldi},~\bfnm{E.~M.}\binits{E.~M.}}
(\byear{2012}).
\btitle{Stochastic blockmodels with a growing number of classes}.
\bjournal{Biometrika}
\bvolume{99}
\bpages{273--284}.
\bid{doi={10.1093/biomet/asr053}, issn={0006-3444}, mr={2931253}}
\bptok{imsref}%
\end{barticle}
%
\endbibitem

\bibitem[\protect\citeauthoryear{Cl{\'e}men{\c{c}}on, Lugosi and
Vayatis}{2008}]{clemenccon2008ranking}
%
\begin{barticle}[mr]
\bauthor{\bsnm{Cl{\'e}men{\c{c}}on},~\bfnm{St{\'e}phan}\binits{S.}},
\bauthor{\bsnm{Lugosi},~\bfnm{G{\'a}bor}\binits{G.}} \AND
\bauthor{\bsnm{Vayatis},~\bfnm{Nicolas}\binits{N.}}
(\byear{2008}).
\btitle{Ranking and empirical minimization of {$U$}-statistics}.
\bjournal{Ann. Statist.}
\bvolume{36}
\bpages{844--874}.
\bid{doi={10.1214/009052607000000910}, issn={0090-5364}, mr={2396817}}
\bptok{imsref}%
\end{barticle}
%
\endbibitem

\bibitem[\protect\citeauthoryear{Diaconis and
Janson}{2008}]{diaconis2007graph}
%
\begin{barticle}[mr]
\bauthor{\bsnm{Diaconis},~\bfnm{Persi}\binits{P.}} \AND
\bauthor{\bsnm{Janson},~\bfnm{Svante}\binits{S.}}
(\byear{2008}).
\btitle{Graph limits and exchangeable random graphs}.
\bjournal{Rend. Mat. Appl. (7)}
\bvolume{28}
\bpages{33--61}.
\bid{issn={1120-7183}, mr={2463439}}
\bptok{imsref}%
\end{barticle}
%
\endbibitem

\bibitem[\protect\citeauthoryear{Fienberg}{2012}]{fienberg2012brief}
%
\begin{barticle}[mr]
\bauthor{\bsnm{Fienberg},~\bfnm{Stephen~E.}\binits{S.~E.}}
(\byear{2012}).
\btitle{A brief history of statistical models for network analysis and open
challenges}.
\bjournal{J. Comput. Graph. Statist.}
\bvolume{21}
\bpages{825--839}.
\bid{doi={10.1080/10618600.2012.738106}, issn={1061-8600}, mr={3005799}}
\bptok{imsref}%
\end{barticle}
%
\endbibitem

\bibitem[\protect\citeauthoryear{Fishkind
et~al.}{2013}]{fishkind2013consistent}
%
\begin{barticle}[mr]
\bauthor{\bsnm{Fishkind},~\bfnm{Donniell~E.}\binits{D.~E.}},
\bauthor{\bsnm{Sussman},~\bfnm{Daniel~L.}\binits{D.~L.}},
\bauthor{\bsnm{Tang},~\bfnm{Minh}\binits{M.}},
\bauthor{\bsnm{Vogelstein},~\bfnm{Joshua~T.}\binits{J.~T.}} \AND
\bauthor{\bsnm{Priebe},~\bfnm{Carey~E.}\binits{C.~E.}}
(\byear{2013}).
\btitle{Consistent adjacency-spectral partitioning for the stochastic block
model when the model parameters are unknown}.
\bjournal{SIAM J. Matrix Anal. Appl.}
\bvolume{34}
\bpages{23--39}.
\bid{doi={10.1137/120875600}, issn={0895-4798}, mr={3032990}}
\bptok{imsref}%
\end{barticle}
%
\endbibitem

\bibitem[\protect\citeauthoryear{Flynn and Perry}{2012}]{flynn2012consistent}
%
\begin{bmisc}[author]
\bauthor{\bsnm{Flynn},~\bfnm{C.~J.}\binits{C.~J.}} \AND
\bauthor{\bsnm{Perry},~\bfnm{P.~O.}\binits{P.~O.}}
(\byear{2012}).
\bhowpublished{Consistent biclustering.
Preprint. Available at \arxivurl{arXiv:1206.6927}.}
\bptok{imsref}%
\end{bmisc}
%
\endbibitem

\bibitem[\protect\citeauthoryear{Fortunato and
Barth{\'{e}}lemy}{2007}]{fortunato2007resolution}
%
\begin{barticle}[pbm]
\bauthor{\bsnm{Fortunato},~\bfnm{Santo}\binits{S.}} \AND
\bauthor{\bsnm{Barth{\'{e}}lemy},~\bfnm{Marc}\binits{M.}}
(\byear{2007}).
\btitle{Resolution limit in community detection}.
\bjournal{Proc. Natl. Acad. Sci. USA}
\bvolume{104}
\bpages{36--41}.
\bid{doi={10.1073/pnas.0605965104}, issn={0027-8424}, pii={0605965104},
pmcid={1765466}, pmid={17190818}}
\bptok{imsref}%
\end{barticle}
%
\endbibitem

\bibitem[\protect\citeauthoryear{Hoeffding}{1963}]{hoeffding1963probability}
%
\begin{barticle}[mr]
\bauthor{\bsnm{Hoeffding},~\bfnm{Wassily}\binits{W.}}
(\byear{1963}).
\btitle{Probability inequalities for sums of bounded random variables}.
\bjournal{J. Amer. Statist. Assoc.}
\bvolume{58}
\bpages{13--30}.
\bid{issn={0162-1459}, mr={0144363}}
\bptok{imsref}%
\end{barticle}
%
\endbibitem

\bibitem[\protect\citeauthoryear{Hoff}{2009}]{hoff2009cmot}
%
\begin{barticle}[author]
\bauthor{\bsnm{Hoff},~\bfnm{P.~D.}\binits{P.~D.}}
(\byear{2009}).
\btitle{Multiplicative latent factor models for description and
prediction of
social networks}.
\bjournal{Computat. Math. Org. Theory}
\bvolume{15}
\bpages{261--272}.
\bptok{imsref}%
\end{barticle}
%
\endbibitem

\bibitem[\protect\citeauthoryear{Hoff, Raftery and
Handcock}{2002}]{hoff2002latent}
%
\begin{barticle}[mr]
\bauthor{\bsnm{Hoff},~\bfnm{Peter~D.}\binits{P.~D.}},
\bauthor{\bsnm{Raftery},~\bfnm{Adrian~E.}\binits{A.~E.}} \AND
\bauthor{\bsnm{Handcock},~\bfnm{Mark~S.}\binits{M.~S.}}
(\byear{2002}).
\btitle{Latent space approaches to social network analysis}.
\bjournal{J. Amer. Statist. Assoc.}
\bvolume{97}
\bpages{1090--1098}.
\bid{doi={10.1198/016214502388618906}, issn={0162-1459}, mr={1951262}}
\bptok{imsref}%
\end{barticle}
%
\endbibitem

\bibitem[\protect\citeauthoryear{Kim and
Leskovec}{2012}]{kim2012multiplicative}
%
\begin{barticle}[mr]
\bauthor{\bsnm{Kim},~\bfnm{Myunghwan}\binits{M.}} \AND
\bauthor{\bsnm{Leskovec},~\bfnm{Jure}\binits{J.}}
(\byear{2012}).
\btitle{Multiplicative attribute graph model of real-world networks}.
\bjournal{Internet Math.}
\bvolume{8}
\bpages{113--160}.
\bid{doi={10.1080/15427951.2012.625257}, issn={1542-7951}, mr={2900491}}
\bptok{imsref}%
\end{barticle}
%
\endbibitem

\bibitem[\protect\citeauthoryear{Miller, Griffiths and
Jordan}{2009}]{miller2009nonparametric}
%
\begin{bincollection}[author]
\bauthor{\bsnm{Miller},~\bfnm{K.~T.}\binits{K.~T.}},
\bauthor{\bsnm{Griffiths},~\bfnm{T.~L.}\binits{T.~L.}} \AND
\bauthor{\bsnm{Jordan},~\bfnm{M.~I.}\binits{M.~I.}}
(\byear{2009}).
\btitle{Nonparametric latent feature models for link prediction}.
In \bbooktitle{Advances in Neural Information Processing Systems 22}
(\beditor{\bfnm{Y.}\binits{Y.}~\bsnm{Bengio}},
\beditor{\bfnm{D.}\binits{D.}~\bsnm{Schuurmans}},
\beditor{\bfnm{J.}\binits{J.}~\bsnm{Lafferty}},
\beditor{\bfnm{C.~K.~I.}\binits{C.~K.~I.}~\bsnm{Williams}} \AND
\beditor{\bfnm{A.}\binits{A.}~\bsnm{Culotta}}, eds.)
\bpages{1276--1284}.
\bpublisher{MIT Press}, \baddress{Cambridge, MA}.
\bptok{imsref}%
\end{bincollection}
%
\endbibitem

\bibitem[\protect\citeauthoryear{Newman}{2006}]{newman2006modularity}
%
\begin{barticle}[author]
\bauthor{\bsnm{Newman},~\bfnm{M.~E.~J.}\binits{M.~E.~J.}}
(\byear{2006}).
\btitle{Modularity and community structure in networks}.
\bjournal{Proc. Natl. Acad. Sci. USA}
\bvolume{103}
\bpages{8577--8582}.
\bptok{imsref}%
\end{barticle}
%
\endbibitem

\bibitem[\protect\citeauthoryear{Rohe, Chatterjee and
Yu}{2011}]{rohe2010spectral}
%
\begin{barticle}[mr]
\bauthor{\bsnm{Rohe},~\bfnm{Karl}\binits{K.}},
\bauthor{\bsnm{Chatterjee},~\bfnm{Sourav}\binits{S.}} \AND
\bauthor{\bsnm{Yu},~\bfnm{Bin}\binits{B.}}
(\byear{2011}).
\btitle{Spectral clustering and the high-dimensional stochastic blockmodel}.
\bjournal{Ann. Statist.}
\bvolume{39}
\bpages{1878--1915}.
\bid{doi={10.1214/11-AOS887}, issn={0090-5364}, mr={2893856}}
\bptok{imsref}%
\end{barticle}
%
\endbibitem

\bibitem[\protect\citeauthoryear{Rohe and Yu}{2012}]{rohe2012co}
%
\begin{bmisc}[author]
\bauthor{\bsnm{Rohe},~\bfnm{K.}\binits{K.}} \AND
\bauthor{\bsnm{Yu},~\bfnm{B.}\binits{B.}}
(\byear{2012}).
\bhowpublished{Co-clustering for directed graphs: The stochastic co-blockmodel
and a
spectral algorithm.
Preprint. Available at \arxivurl{arXiv:1204.2296}.}
\bptok{imsref}%
\end{bmisc}
%
\endbibitem

\bibitem[\protect\citeauthoryear{Schneider}{1993}]{schneider1993convex}
%
\begin{bbook}[mr]
\bauthor{\bsnm{Schneider},~\bfnm{Rolf}\binits{R.}}
(\byear{1993}).
\btitle{Convex Bodies: The {B}runn--{M}inkowski Theory}.
\bseries{Encyclopedia of Mathematics and Its Applications}
\bvolume{44}.
\bpublisher{Cambridge Univ. Press}, \blocation{Cambridge}.
\bid{doi={10.1017/CBO9780511526282}, mr={1216521}}
\bptok{imsref}%
\end{bbook}
%
\endbibitem

\bibitem[\protect\citeauthoryear{White}{1982}]{white1982maximum}
%
\begin{barticle}[mr]
\bauthor{\bsnm{White},~\bfnm{Halbert}\binits{H.}}
(\byear{1982}).
\btitle{Maximum likelihood estimation of misspecified models}.
\bjournal{Econometrica}
\bvolume{50}
\bpages{1--25}.
\bid{doi={10.2307/1912526}, issn={0012-9682}, mr={0640163}}
\bptok{imsref}%
\end{barticle}
%
\endbibitem

\bibitem[\protect\citeauthoryear{Zhao, Levina and
Zhu}{2011}]{zhao2011community}
%
\begin{barticle}[pbm]
\bauthor{\bsnm{Zhao},~\bfnm{Yunpeng}\binits{Y.}},
\bauthor{\bsnm{Levina},~\bfnm{Elizaveta}\binits{E.}} \AND
\bauthor{\bsnm{Zhu},~\bfnm{Ji}\binits{J.}}
(\byear{2011}).
\btitle{Community extraction for social networks}.
\bjournal{Proc. Natl. Acad. Sci. USA}
\bvolume{108}
\bpages{7321--7326}.
\bid{doi={10.1073/pnas.1006642108}, issn={1091-6490}, pii={1006642108},
pmcid={3088589}, pmid={21502538}}
\bptok{imsref}%
\end{barticle}
%
\endbibitem

\bibitem[\protect\citeauthoryear{Zhao, Levina and
Zhu}{2012}]{zhao2013consistency}
%
\begin{barticle}[mr]
\bauthor{\bsnm{Zhao},~\bfnm{Yunpeng}\binits{Y.}},
\bauthor{\bsnm{Levina},~\bfnm{Elizaveta}\binits{E.}} \AND
\bauthor{\bsnm{Zhu},~\bfnm{Ji}\binits{J.}}
(\byear{2012}).
\btitle{Consistency of community detection in networks under degree-corrected
stochastic block models}.
\bjournal{Ann. Statist.}
\bvolume{40}
\bpages{2266--2292}.
\bid{doi={10.1214/12-AOS1036}, issn={0090-5364}, mr={3059083}}
\bptok{imsref}%
\end{barticle}
%
\endbibitem

\end{thebibliography}
\end{document}